\documentclass[12pt,reqno]{amsart}

\usepackage{amssymb}
\usepackage[T1]{fontenc}
\usepackage[latin1]{inputenc} 
\usepackage{dsfont, mathrsfs}
\usepackage{stmaryrd}
\usepackage{amsthm}
\usepackage{graphicx}
\usepackage{color, xcolor, soul}
\usepackage{xspace}
\usepackage{yhmath, epigraph}
\usepackage{caption, subcaption}
\usepackage{relsize}
\usepackage{shuffle}
\usepackage{float, adjustbox} 
\usepackage{mleftright}
\usepackage{pdfpages}
\usepackage{pax} 
\usepackage[
        colorlinks,
        hyperfootnotes=true,
        pagebackref, 
        hyperindex,
        linkcolor=blue]{hyperref}
\usepackage{geometry}
\usepackage{framed} 
\colorlet{shadecolor}{gray!40}
\usepackage{bbm}
\usepackage[bbgreekl]{mathbbol}
\DeclareSymbolFontAlphabet{\mathbb}{AMSb}
\DeclareSymbolFontAlphabet{\mathbbl}{bbold}
\usepackage{tikz}
\usetikzlibrary{positioning}
\usetikzlibrary{calc}
\usetikzlibrary{arrows.meta}
\usetikzlibrary{arrows}
\usepackage{pgfplots}
\pgfplotsset{compat=1.18}
\makeatletter 
\newcommand\HREF[2]{\hyper@linkurl{#2}{#1}}
\makeatother  
\usepackage{verbatim}
\DeclareFontFamily{U}{mathx}{\hyphenchar\font45}
\DeclareFontShape{U}{mathx}{m}{n}{
      <5> <6> <7> <8> <9> <10>
      <10.95> <12> <14.4> <17.28> <20.74> <24.88>
      mathx10
      }{}
\DeclareSymbolFont{mathx}{U}{mathx}{m}{n}
\DeclareFontSubstitution{U}{mathx}{m}{n}
\DeclareMathAccent{\widecheck}{0}{mathx}{"71}
\DeclareMathAccent{\wideparen}{0}{mathx}{"75}

\newcommand{\Zindex}[1]{}  
\usepackage{listings} 
\usepackage{marginnote}
\allowdisplaybreaks

\theoremstyle{plain}
\newtheorem{theorem}{Theorem}[section]
\newtheorem{lemma}[theorem]{Lemma}

\newtheorem{proposition}[theorem]{Proposition}

\theoremstyle{definition}
\newtheorem{definition}[theorem]{Definition}

\newtheorem{question}[theorem]{Question}

\theoremstyle{remark}
\newtheorem{remark}[theorem]{Remark}

\newcommand{\Cc}{\mathbb{C}}

\newcommand{\Ee}{\mathbb{E}}
\newcommand{\Ff}{\mathbb{F}}

\newcommand{\Mm}{\mathbb{M}}
\newcommand{\Nn}{\mathbb{N}}
\newcommand{\Pp}{\mathbb{P}}

\newcommand{\Rr}{\mathbb{R}}

\newcommand{\Uu}{\mathbb{U}}

\newcommand{\Yy}{\mathbb{Y}}
\newcommand{\Zz}{\mathbb{Z}}

\newcommand{\Un}{\mathds{1}}

\newcommand{\Ae}{\mathcal{A}}
\newcommand{\Be}{\mathcal{B}}

\newcommand{\Fe}{\mathcal{F}}

\newcommand{\Le}{\mathcal{L}}

\newcommand{\Pe}{\mathcal{P}}

\newcommand{\Ue}{\mathcal{U}}

\newcommand{\Ze}{\mathcal{Z}}

\newcommand{\Beb}{{{\boldsymbol{\mathcal{B}}}}}

\newcommand{\Deb}{{{\boldsymbol{\mathcal{D}}}}}
\newcommand{\Geb}{{{\boldsymbol{\mathcal{G}}}}}
\newcommand{\Heb}{{{\boldsymbol{\mathcal{H}}}}}

\newcommand{\Neb}{{{\boldsymbol{\mathcal{N}}}}}
\newcommand{\Peb}{{{\boldsymbol{\mathcal{P}}}}}

\newcommand{\Ns}{{\mathscr{N}}}

\newcommand{\Us}{{\mathscr{U}}}

\newcommand{\Sg}{{\mathfrak{S}}}

\newcommand{\deltab}{{\boldsymbol{\delta}}}
\newcommand{\gammab}{{\boldsymbol{\gamma}}}

\newcommand{\lambdab}{{\boldsymbol{\lambda}}}

\newcommand{\pib}{{\boldsymbol{\pi}}}

\newcommand{\sigmab}{{\boldsymbol{\sigma}}}

\newcommand{\zetab}{{\boldsymbol{\zeta}}}

\newcommand{\Gammab}{{\boldsymbol{\Gamma}}}

\newcommand{\Cb}{{\boldsymbol{C}}}

\newcommand{\ab}{{\boldsymbol{a}}}

\newcommand{\cb}{{\boldsymbol{c}}}

\newcommand{\vb}{{\boldsymbol{v}}}

\newcommand{\xb}{{\boldsymbol{x}}}

\def\ee{ { \mathbbm{e} } }

\newcommand{\ensemble}[1]{ {\left\lbrace #1 \right\rbrace } } 
\newcommand{\prth}[1]{\!\left( #1 \right) }
\newcommand{\crochet}[1]{\!\left[ #1 \right] }  
\newcommand{\intcrochet}[1]{\llbracket #1 \rrbracket} 
\newcommand{\abs}[1]{{\left| #1 \right|}}  
\newcommand{\norm}[1]{\left| \! \left| #1 \right| \! \right|}

\newcommand{\pleth}[1]{{\left[ #1 \right] }}

\newcommand{\Esp}[1]{ \Ee  \prth{ #1 } }  
\newcommand{\Prob}[1]{ \Pp \prth{ #1 } } 
 
\def\inv{^{-1}}

\newcommand{\cvlaw}[2]{\xrightarrow[#1 \, \rightarrow \, #2]{\Le}} 
\newcommand{\tendvers}[2]{\xrightarrow[#1 \, \rightarrow \, #2]{}} 
\newcommand{\cvmodp}[2]{\xrightarrow[#1 \, \rightarrow \, #2]{\operatorname{mod-P}}}

\newcommand{\equivalent}[1]{ {\underset{#1 }{\sim} } }

\def\divise{ \vert }
\def\mobius{\mu}

\newcommand{\Unens}[1]{ \Un_{ \ensemble{#1} } }
\newcommand{\pe}[1]{{\left[ #1 \right]}}
\newcommand{\pf}[1]{{\left\{ #1 \right\}}}

\def\eqlaw{\stackrel{\Le}{=}}

\def\card{{ \operatorname{card} }}
\def\str{{ \operatorname{str} }}
\def\Cent{{ \operatorname{Cent} }}

\def\Exp{ \operatorname{Exp} }
\def\Ber{ \operatorname{Ber} }
\def\Geom{ \operatorname{Geom} }
\def\Poisson{ \operatorname{Po} } 

\def\NegBin{\operatorname{NegBin}}
\def\Bin{\operatorname{Bin}}

\def\geq{\geqslant}
\def\leq{\leqslant}

\def\Re{\mathfrak{Re}}

\let\oldforall\forall
\def\forall{\oldforall\,} 

\let\oldexists\exists
\def\exists{\oldexists\,}

\newcommand{\emailhref}[1]{ \email{\href{mailto:#1}{#1}} }

\def\TAMS{Trans. Amer. Math. Soc. }

\title[Splitting phenomenon and universality of the $ \Gamma $ factor ]{Splitting phenomenon in the Sathe-Selberg theorem, mod-Poisson convergence with auxiliary randomisation and universality of the Gamma factor}
\author[Y. Barhoumi-Andr\'eani]{Yacine Barhoumi-Andr\'eani}
\emailhref{yacine.barhoumi@math.bas.bg}
\address{Bulgarian Academy of Sciences,       
Institute of Mathematics and Informatics, 
Department of Algebra and Logic, 
Acad. Georgi Bonchev Str., Block 8, 1113 Sofia (Bulgaria).}
\date{\today}


\begin{document}
\begin{abstract}
We consider several sequences of random variables whose Fourier-Laplace transforms present the same type of \textit{splitting phenomenon} when suitably rescaled by the Fourier-Laplace transform of a Poisson-distributed random variable (mod-Poisson convergence). Addressing a question raised by Kowalski-Nikeghbali, we explain the appearance of a universal term, the \textit{Gamma factor}, by a common feature of each model, the existence of an auxiliary randomisation that reveals an independence structure. This universal Gamma factor is a consequence of the probabilistic (exponential) fluctuations of this subordinated random variable.

The class of examples that belong to this framework includes: random uniform permutations, random Gauss integers, random polynomials over a finite field, random matrices with values in a finite field, random partitions and the classical Sathe-Selberg theorems in probabilistic number theory. For this last example, the randomisation that triggers the Gamma factor defines a new probability distribution, the ``delta-Zeta distribution'' which is the analogue of the geometric distribution in the permutation setting.
%
%
\end{abstract}
\maketitle

\setcounter{tocdepth}{1}
\tableofcontents



\section{Introduction}\label{Sec:Intro}

\subsection{Motivations}

\subsubsection{\textbf{The splitting phenomenon}}

The Sathe-Selberg theorems in probabilistic number theory concern the \textit{second-order fluctuations} of the arithmetic random variables given by the number of prime divisors $ \omega(U_n) $ of a random uniform integer $ U_n \in \ensemble{1, 2, \dots, n} $, or its number of prime divisors counted with multiplicities $ \Omega(U_n) $. 

Let $ \Pe $ denote the set of prime numbers. For $ k \geq 1 $ denote by $ v_p(k) $ the $p$-adic valuation of $ k $, namely the exponent of $p$ in its prime decomposition:
\begin{align}\label{Thm:PrimeDecomposition}
k = \prod_{p \in \Pe } p^{ v_p(k) }
\end{align}

The number of prime divisors of $k$ and the number of prime divisors counted with multiplicities are respectively defined by
\begin{align*}
\omega(k) :=  \sum_{p \in \Pe } \Unens{ v_p(k) \geq 1 }, 
\qquad
\Omega(k) :=  \sum_{p \in \Pe } v_p(k)
\end{align*}
We now consider the ``arithmetic'' random variables $ \omega(U_n) $ and $ \Omega(U_n) $ for $ U_n \sim \Ue\prth{ \ensemble{1, 2, \dots, n} } $. 
The classical Erd\"os-Kac theorems (see \cite{ErdosKac, ErdosKac2, RenyiTuran}) can be stated as
\begin{align}\label{Thm:ErdosKac}
\sup_{x \in \Rr} \abs{ \Prob{ \frac{ X_n - \log\log n }{ \sqrt{ \log\log n }  }  \leq x } - \int_{-\infty }^x e^{- u^2 / 2 } \frac{ du }{ \sqrt{2 \pi } } } 
\tendvers{n}{ + \infty } 0,  
\qquad
X_n \in \ensemble{\omega(U_n), \Omega(U_n)} 
\end{align}

A more refined version of this last theorem shows in fact that each random variable $ X_n $ can be approximated by a Poisson random variable $ \Peb(\log\log n) \sim \Poisson(\log\log n) $:
\begin{align}\label{Thm:ErdosKacPoisson}
\sup_{x \in \Rr} \abs{ \Prob{ \frac{ X_n - \log\log n }{ \sqrt{ \log\log n }  }  \leq x } -  \Prob{ \Peb(\log\log n) \leq x} } \tendvers{n}{ + \infty } 0 
\end{align}

This Central Limit Theorem/Poisson approximation can be regarded as a computation of \textit{first order} fluctuations of $ \omega(U_n)$ and $ \Omega(U_n) $. In particular, \eqref{Thm:ErdosKac} is a convergence in distribution that is equivalent by L\'evy injectivity to a convergence of the Fourier or Laplace transform.


The Sathe-Selberg theorems (see \cite{Sathe, SelbergSathe}) are a refinement of this last convergence and can be regarded as \textit{second order} fluctuations of these last quantities. They read, locally uniformly in $ x \in \Cc $
%
%
%
\begin{align}\label{Thm:SatheSelberg}
\frac{ \Esp{ x^{  X_n } } }{  \Esp{ x^{ \Peb( \gamma_n ) } } } \tendvers{n}{ + \infty } \Phi_X\prth{ x }
\end{align}
where $ \gamma_n = \log\log n + O(1) $ and $ \Phi_X : \Cc \to \Rr_+ $ is a function depending on the sequence $ (X_n)_n $. Writing $ \Phi_\omega $ for $ X_n = \omega(U_n) $ and $ \Phi_\Omega $ for $ X_n = \Omega(U_n) $, one gets, for all $ x \in (0, 1) $
\begin{align}\label{Def:FonctionsModPoissonSatheSelberg}
\begin{aligned}
\Phi_\omega(x) & =  \frac{ 1}{ \Gamma(x) } \prod_{p \in \Pe} \prth{ 1 + \frac{x-1}{p} } e^{ - \frac{x-1}{p} }    \\
\Phi_\Omega(x) & =  \frac{1}{ \Gamma(x) } \prod_{p \in \Pe} \frac{1 - \frac{1}{p} }{ 1 - \frac{x }{p} } e^{ - \frac{x-1}{p} } 
\end{aligned}
\end{align}

Note that these functions depend on the choice of the constant $ O(1) $ in $ \gamma_n $, but all the considered functions are equal, up to a multiplication by $ e^{ \kappa (x-1) } = \Esp{ x^{\Peb(\kappa) } } $. The type of renormalisation that occurs in \eqref{Thm:SatheSelberg} was hence called \textit{mod-Poisson} convergence by Kowalski and Nikeghbali \cite{KowalskiNikeghbali2} who introduced it following a similar study in the Gaussian setting \cite{JacodAl} (see also \cite{BarbourKowalskiNikeghbali, Hwang96PGD, Hwang98PGD, Hwang98rate}).

The limiting functions that occur in \eqref{Thm:SatheSelberg} remember certain information from the Erd\"os-Kac theorem \eqref{Thm:ErdosKac}: consider the \textit{independent models} 
\begin{align}\label{Def:IndependentModel}
\widehat{\omega}_n  := \sum_{ p \in \Pe\! , \, p \leq n } \Beb\prth{ \frac{1}{p} }, 
\qquad 
\widehat{\Omega}_n  := \sum_{ p \in \Pe\! , \, p \leq n } \Geb\prth{ \frac{1}{p} }
\end{align}
where the $ (\Beb(1/p))_p $ and the $ (\Geb(1/p))_p $ are independent random variables with Bernoulli and geometric distribution respectively.

In a certain sense, these last random variables are good surrogates of their original counterparts as they succeed in explaining heuristically the Erd\"os-Kac theorem: using $ \sum_{p \leq n} p\inv = \log\log n + O(1) $ and applying the classical CLT for sums of independent random variables, one gets \eqref{Thm:ErdosKac}. These independent models support moreover the Poisson approximation \eqref{Thm:ErdosKacPoisson} using the classical Poisson approximation for sums of independent Bernoulli random variables. The heuristic replacement of the dependent Bernoulli (resp. Geometric) random variables by their independent counterparts is thus a full success at the level of the fluctuations, and this, despite having \textit{totally dependent} random variables constructed out of the same exact source of randomness $ U_n $.

These independent models nevertheless fail to completely explain the form of the limiting functions \eqref{Def:FonctionsModPoissonSatheSelberg} since 
\begin{align}\label{Thm:CvModPoissonModelesIndependants}
\begin{aligned}
\frac{ \Esp{ x^{ \, \widehat{\omega}_n } } }{  \Esp{ x^{ \Peb(\gamma_n) } } } 
& \tendvers{n}{ + \infty } \Phi_{\widehat{\omega}}(x) := \prod_{p \in \Pe} \prth{ 1 + \frac{x-1}{p} } e^{ - \frac{x-1}{p} } \\
\frac{ \Esp{ x^{ \, \widehat{\Omega}_n } } }{  \Esp{ x^{ \Peb(\gamma_n) } } } 
& \tendvers{n}{ + \infty } \Phi_{\widehat{\Omega}}(x) := \prod_{p \in \Pe} \frac{1 - \frac{1}{p} }{ 1 - \frac{x }{p} } e^{ - \frac{x-1}{p} }  
\end{aligned}
\end{align}

Compared to \eqref{Thm:CvModPoissonModelesIndependants}, the key point to notice in \eqref{Def:FonctionsModPoissonSatheSelberg} is thus the appearance of the common \textit{Gamma factor} $ 1/\Gamma(x) $ as a correction to the independent approximation:
\begin{align}\label{Eq:SplittingPoisson}
\Phi_\omega = \frac{1}{\Gamma}\times \Phi_{\widehat{\omega}}, 
\qquad
\Phi_\Omega = \frac{1}{\Gamma}\times \Phi_{\widehat{\Omega}}
\end{align}


%

\medskip
\subsubsection{\textbf{Universality, 1}}

The persistence of such a factor in all these models is reminiscent of the \textit{moments conjecture} in Probabilistic Number Theory: the random variable $ \Ze_t := \log\abs{ \zeta\prth{ \frac{1}{2} + i t U } } $ with $ U $ a random variable uniform in $ \crochet{0, 1} $ satisfies a Gaussian CLT originally due to Selberg (see e.g. \cite{TaoBlog})
\begin{align}\label{Thm:SelbergCLT}
\frac{ \Esp{ e^{- \lambda \Ze_t / \sqrt{ \frac{1}{2} \log\log t } } } }{ \Esp{e^{-\lambda G} } }
\tendvers{t}{ + \infty } 1 
\end{align}
with $ G \sim \Ns(0, 1) $, and the second order fluctuations are conjecturally given by the following \textit{mod-Gaussian convergence} \cite{JacodAl}:
\begin{align}\label{Thm:KeatingSnaithMomentsConjecture}
\frac{ \Esp{ e^{- \lambda \Ze_t } } }{ \Esp{ e^{ -\lambda \sqrt{ \frac{1}{2} \log\log t } \, G } } }
\tendvers{t}{ + \infty }  
\Phi_\zeta(\lambda) := \Phi_{\widehat{\zeta}}(\lambda) \Phi_{\mathrm{RMT}}(\lambda)   
\end{align}

This is the \textit{moments conjecture} \cite{KeatingSnaith}. The factor $ \Phi_{\widehat{\zeta}} $ is particular to the considered model and comes from an independent model where each term in the prime product that defines the Riemann $ \zeta $ function inside $ \ensemble{\Re > 1} $ has been replaced by an independent counterpart (the exact same previous heuristic that led to $ \widehat{\omega}_n $ and $ \widehat{\Omega}_n $). 
The factor $ \Phi_{\mathrm{RMT}} $  is supposedly universal and was conjectured in relation with Random Matrix Theory ; it comes from the mod-Gaussian convergence of the characteristic polynomial on the unit circle of a random Haar-distributed unitary matrix \cite{KeatingSnaith2, BHNY}.

The mod-Gaussian convergence \eqref{Thm:KeatingSnaithMomentsConjecture} is one instance of a splitting phenomenon: the limiting function is a product of two terms, one of which is universal in the sense that it appears in other (a priori unrelated) models. It is in fact quite common in the study of functionals of random unitary matrices, see e.g. the work \cite{BarhoumiCUErevisited} by the author where several other such models display the same type of behaviour, or the explanation by Nikeghbali-Yor \cite{NikeghbaliYorBarnes} of a similar splitting phenomenon in the function $ \Phi_{\mathrm{RMT}} : z \mapsto G(1 + z)^2 \times G(1 + 2z)\inv $ using a decomposition of the random variables at stake. The phenomenon present in the Sathe-Selberg theorem \eqref{Def:FonctionsModPoissonSatheSelberg} could be of the same type, since the limiting mod-Poisson function splits into a product involving the Gamma factor and the independent model factor. But to really motivate the terminology of ``universal factor'', it remains to explicit a zoo of models that share this univerality feature.

\subsubsection{\textbf{Universality, 2}}

The \textit{splitting phenomenon} occurring in \eqref{Eq:SplittingPoisson} is also present in other discrete models that converge in the mod-Poisson sense. 
For instance, let $ P_n $ denote a random polynomial uniformly distributed in $ \ensemble{ P \in \Ff_q[X]_m :  \deg P = n  } $ where $ q $ is a certain power of a prime number and $ \Ff_q[X]_m $ designates the monic such polynomials with coefficients in $ \Ff_q $. Let $ \Pe\prth{  \Ff_q[X] } $ denote the set of irreducible monic polynomials of $ \Ff_q[X] $ and for $ Q \in \Ff_q[X]_m $, define the number of irreducible divisors of $ Q $ counted without multiplicity by 
\begin{align*}
\omega_q(Q) := \sum_{ \pi \in \Pe(\Ff_q[X]) } \Unens{ \pi \divise Q }
\end{align*}

One could define in the same vein $ \Omega_q(Q) $ using the prime decomposition analogous to \eqref{Thm:PrimeDecomposition} in the function field setting. A classical result then shows that (see e.g. \cite[thm. 6.1]{KowalskiNikeghbali2})
\begin{align*}
\frac{ \Esp{ x^{ \, \omega_q(P_n) } } }{ \Esp{  x^{ \, \Peb(\gamma_{n, q}) } } }  
\tendvers{n }{ + \infty }  
\Phi_{\omega_q}\prth{ x }
\end{align*}
with
\begin{align}\label{Def:FonctionModPoissonPolynomes}
\Phi_{\omega_q}(x) = \frac{1}{\Gamma(x) } \times  \prod_{ \pi \in \Pe(\Ff_q[X]) } \prth{ 1 + \frac{x-1}{ \abs{\pi}_q } } e^{ - \frac{x-1}{ \abs{\pi}_q } }
\end{align}
where $ \abs{\pi}_q = q^{ \deg(\pi) } $ and $ \gamma_{n,  q} = \log n + O(1) $. Here again, one sees that an independent model made of independent Bernoulli random variables would give a function $ \Phi_{\widehat{\omega}_q} := \Phi_{\omega_q} \times \Gamma $.  

The list of models that incorporate such a \textit{Gamma factor} in a mod-Poisson convergence includes:
\begin{itemize}

\item the total number of cycles $ C(\sigmab_{\! n}) $ of a random uniform permutation $ \sigmab_{\! n} \in \Sg_n $ the symmetric group acting on $n$ elements, and more precisely~:
\begin{align}\label{Thm:CvModPoissonCyclesUnifPerm}
\frac{ \Esp{ x^{ \, C(\sigmab_{\!n}) } } }{  \Esp{ x^{ \Peb(\ln n) } } } 
& \tendvers{n}{ + \infty } \Phi_C(x) := \frac{1}{\Gamma(x)}
\end{align}

\medskip
\item the number of prime divisors of a random uniform Gauss integer (the analogue of the Sathe-Selberg theorems in $ \Zz\crochet{i} $),

\medskip
\item the total number of Jordan blocks $ C(\pib_{n, q}) $ of a random uniform matrix of $ \pib_{n, q} \in GL_n(\Ff_q) $, 

\medskip
\item \&c. 

\end{itemize}

A comprehensive and non exhaustive table of such models is recapitulated in \S~\ref{SubSec:Conclusion:Summary}.

\medskip
\subsubsection{\textbf{Universality, 3}}

One of the main motivations of E. Kowalski and A. Nikeghbali for introducing mod-Poisson convergence in \cite{KowalskiNikeghbali2} after having introduced mod-Gaussian convergence with J. Jacod in \cite{JacodAl} was the existence of this variety of instances of splitting phenomena in probability theory. This naturally led them to ask for an even more general type of universality when concerned with splitting of limiting functions (and later beyond the Gaussian and Poisson case with the concept of mod-$ \phi $ convergence \cite{DelbaenAl, FerayMeliotNikeghbali1, FerayMeliotNikeghbali2}). The \textit{Kowalski-Nikeghbali philosophy} is the search for models/random variables that share this property with the idea that a general paradigm should be at stake to explain such a systematic occurrence. It is a refinement of the \textit{Kac-Sarnak philosophy} \cite{KatzSarnak} that studies models related to finite fields (more precisely Dirichlet characters mod $ t^N $ on $ \Ff_q[t] $ when $ q \to +\infty $ with fixed $N$) and the \textit{Keating-Snaith philosophy} \cite{KeatingSnaith, KeatingSnaith2} that studies the characteristic polynomial of a random Haar-distributed random unitary matrix as a toy-model for the Riemann $ \zeta $ function, both philosophies related amongst others with splitting phenomena in mod-Gaussian convergence.

So far, two paradigms of explanation for the splitting phenomenon were provided in the mod-Gaussian setting of $ \log\abs{\zeta(\frac{1}{2} + iT U)} $: the \textit{hybrid product paradigm} developed in \cite{GonekHughesKeating} and the \textit{multiple Dirichlet series paradigm} developed in \cite{DiaconuGoldfeldHoffstein}. In the mod-Poisson framework, the original paradigm in \cite{KowalskiNikeghbali2} treats the case of $ \omega_q(P_n) $ using a sheaf-theoretic structure on $ \Ff_q[t] $ and an approximation of random polynomials by quadratfrei polynomials, leading to a Galois group equal to the symmetric group $ \Sg_n $ and triggering so in a natural way $ C(\sigmab_{\!n}) $~; another paradigm was proposed for $ \omega(U_n) $ in \cite{BarhoumiModOmega} using the conditioning of a random subset of primes to be present a.s. in the prime factor decomposition, the imposed random sequence of primes involving a \textit{paintbox process} naturally linked to random uniform permutations and to the limit \eqref{Thm:CvModPoissonCyclesUnifPerm}. 


\medskip 
\subsection{Goal and main result}

The previous results lead to the natural questions:

\begin{question}[Kowalski-Nikeghbali \cite{KowalskiNikeghbali2}]\label{Q:KN}
Can one find an interesting paradigm that would explain the universal appearance of the Gamma factor in all models converging in the mod-Poisson sense with such a splitting feature~?
\end{question}


\begin{question}[Key Question]\label{Q:KN++}
Can one find a paradigm that would explain the universal appearance of all universal factors in all models converging in the mod-$ \phi $ sense with such a splitting feature~? I.e. such a paradigm would be the same in the mod-Gaussian and mod-Poisson case.
\end{question}

\medskip

The goal of this article is to bring a particular answer to question~\ref{Q:KN}, i.e. to identify a common mechanism that systematically explains the emergence of the Gamma factor in mod-Poisson convergence. Our answer will be of probabilistic nature and will involve a hidden characteristic of all considered models: the existence of a particular randomisation that turns them into an \textit{independent model} different from \eqref{Def:IndependentModel}. These different randomisations will depend on a parameter and will satisfy a CLT when this parameter tends to a certain value, with a limit in law given by a random variable with exponential distribution. The key idea is that these limiting exponential fluctuations, whose Mellin transform is the Gamma function, will be responsible for the appearance and the universality of the Gamma factor.


\medskip


While we would love to also give here an answer to question~\ref{Q:KN++}, honesty must make us admit that this goal is still far from being achieved, and that the author's attempts to extend it to mod-Gaussian models have not proved successful. We will nevertheless come back to it in a near future with a different paradigm that shall provide such an answer.

\medskip 
\subsection{Organisation of the article}

The plan of the paper is as follows:
\begin{itemize}

\item in \S~\ref{Sec:Notations}, we introduce some notations and conventions used throughout the text; 

\medskip
\item in \S~\ref{Sec:Sn}, we develop the case of the symmetric group that will give one first example of the randomisation paradigm with all the relevant key concepts that are necessary to trigger the Gamma factor (CIP, randomisation with exponential fluctuations, etc.), 

\medskip
\item in \S~\ref{Sec:Integers}, we develop the key example of $ \omega(U_n) $ and introduce the relevant randomisation at stake in this framework, the \textit{Delta-Zeta distribution}, 

\medskip
\item in \S~\ref{Sec:GammaUniversality}, we build on the two previous sections to describe a zoo of functionals whose evaluation in a particular uniformly distributed random variable will give a Gamma factor multiplied by an independent model factor in a mod-Gaussian convergence,

\medskip
\item in \S~\ref{Sec:Conclusion}, we summarise the results in a table in \S~\ref{SubSec:Conclusion:Summary} and discuss some perspectives of interest,

\medskip
\item in Annex~\ref{Sec:DeltaZeta}, we describe some properties of the Delta-Zeta distribution that allow to compare it to the Geometric distribution in a more systematic way.
\end{itemize}


\medskip
\section{Notations and conventions}\label{Sec:Notations}

\subsection{General notations}

We gather here some notations used throughout the text. We will use the multi-index notation $ \xb^\ab := \prod_{k \geq 1} x_k^{a_k}  $. The ``tensor/direct sum'' $ A \oplus B $ designates the diagonal block concatenation of the matrices $ A $ and $B$. The Abel summation designates the equality for summable sequences $ (u_k)_k $ and $ (v_k)_k $ such that $ u_0 := 0 $
\begin{eqnarray*}
\sum_{k \geq 1} u_k (v_k - v_{k + 1}) = \sum_{k \geq 1}  (u_k - u_{k - 1}) v_k
\end{eqnarray*}

We define $ x^{ \uparrow k } := x(x + 1) \dots (x + k - 1) $ and $ x^{\downarrow k} := x(x - 1) \cdots (x - k + 1) $ if $ k \in \Nn $.

\subsubsection{\textbf{Probability}}

We fix a probability space $ (\Omega, \Fe, \Pp) $ where all encountered random variables will be considered. 
The uniform distribution on a set $ X $ will be denoted by $ \Us(X) $, 
the Poisson distribution of parameter $ \gamma > 0 $ by $ \Poisson(\gamma) $, 
the geometric distribution of parameter $ t \in (0, 1) $ by $ \Geom(t) $, 
the exponential distribution of parameter $ t > 0 $ by $ \Exp(t) $,
the Bernoulli distribution of parameter $ q \in (0, 1) $ by $ \Ber_{\{0, 1\}}(q) $ (or $\Ber_{\{\pm 1\}}(q) $ if the random variable has values in $ \{\pm 1 \} $), 
the Binomial distribution of parameter $ n \in \Nn^* $ and $ p \in [0, 1] $ by $ \Bin(n, p) $ 
and the negative binomial distribution of parameters $ s > 0 $ and $ t \in (0, 1) $ by $ \NegBin(s, t) $. 
The fact that a random variable $ Z $ has a certain distribution $ \Ze  $ in the will be denoted by $ Z \sim \Ze $. We recall in particular that for all $ k \in \Nn $ and $ \theta \in \Rr $ 
\begin{align*}
\Peb(\gamma) \sim \Poisson(\gamma) &
\quad\Longrightarrow\quad
\Pp(\Peb(\gamma) = k) = e^{- \gamma } \frac{\gamma^k}{k!} , 
\quad
\Esp{e^{i \theta \Peb(\gamma) } } = \exp\prth{ \gamma ( e^{i\theta } - 1 ) }  \\
\Geb(t) \sim \Geom(t) &
\quad\Longrightarrow\quad
\Pp(\Geb(t) = k) = (1 - t) t^k , 
\quad
\Esp{e^{i \theta \Geb(t) } } = \frac{1 - t}{1 - t e^{i\theta} }  \\
\Neb\!\Beb(t, s) \sim \NegBin(t, s) &
\quad\Longrightarrow\quad
\Pp(\Neb\!\Beb(t, s)  = k) = (1 - t)^s \frac{t^k \, s^{\uparrow k}}{k!}  , \\
& \hspace{+7cm}
\Esp{e^{i \theta \Neb\!\Beb(t, s)   } } = \prth{ \frac{1 - t}{1 - t e^{i\theta} }}^{\!\! s}  \\
\Beb(p) \sim \Ber_{\{0, 1\}}(p) &
\quad\Longrightarrow\quad
\Pp(\Beb(p) = k) = p\Unens{k = 1} + (1 - p) \Unens{k = 0} , \\
& \hspace{+7cm}
\Esp{e^{i \theta \Beb(p) } } = 1 + (p - 1) e^{i\theta}  \\
\Beb(n, p) \sim \Bin(n, p) &
\quad\Longrightarrow\quad
\Pp(\Beb(n,p) = k) = p^k (1 - p)^{n - k} {n \choose k} , \\
& \hspace{+7cm}
\Esp{e^{i \theta \Beb(n, p) } } = \prth{ 1 + (p - 1) e^{i\theta}  }^n
\end{align*}

\subsubsection{\textbf{Properties of the Geometric distribution}}

We recall that if $ \ee \sim \Exp(1) $, $ \Prob{\ee \geq x} = e^{-x} $ for all $ x \geq 0 $ and
\begin{align}\label{Eq:TransfoMellinExp}
\boxed{\Esp{ \ee^{x - 1} } = \Gamma(x)}
\end{align}

The Geometric distribution enjoys the following classical properties whose proof is left as an easy exercise~:
\begin{enumerate}

\item $ \Geb(t) \eqlaw \pe{ \frac{\ee }{ \log(t\inv) }  } $ with $ \ee \sim \Exp(1) $, $ t < 1 $,
%
%
%

\item $  \log(t\inv) \, \Geb(t) \cvlaw{t}{ 1^- } \ee  $,

\item $ \Geb(t) \eqlaw \Peb\prth{ \frac{t}{1 - t} \, \ee } $ with 
$ \ee \sim \Exp(1) $, $ \Peb(\cdot) \sim \Poisson(\cdot) $ and where the randomisation by $ \ee $ has to be understood as the independent randomisation of the parameter of the underlying random variable, 

%
\end{enumerate}

\subsubsection{\textbf{Sets}}

The set $ \ensemble{1, 2, \dots, n } $ will be denoted by $ \intcrochet{1, n} $, the set of prime numbers will be denoted by $ \Pe $, the set of matrices over a finite field $ \Ff_q $ (for a prime power $q$) by $ M_n(\Ff_q) $, the set of invertible such matrices by $ GL_n(\Ff_q) $ and the symmetric group acting on $n$ letters will be denoted by $ \Sg_n $. The set of polynomials with coefficients over a finite field will be denoted by $ \Ff_q[X] $, the set of such monic polynomials by $ \Ff_q[X]_m $, the set of irreducible (or prime) monic such polynomials by $ \Phi_q $. Last, we will also set $ \Phi_q^+ := \Phi_q\setminus\ensemble{X } $.

\subsubsection{\textbf{Functions and constants}}

The Mobius function of $ \Nn $ is denoted by $ \mobius $ and is defined by $ \mu(n) := (-1)^{\omega(n) } \Unens{ \Omega(n) = \omega(n) } $. The integer part of $ x \in \Rr $ will be denoted by $ \pe{x} $. We set $ H_n := \sum_{1 \leq k \leq n} k\inv = \ln(n) + \gamma + o(1) $ and $ H_n^{(\Pe)} := \sum_{ p \in \Pe, p \leq n } p\inv = \ln\ln n + \gamma_\Pe + o(1) $ where $ \gamma $ and $ \gamma_\Pe $ are the Euler and the prime Euler constants.

\subsubsection{\textbf{Partitions}}

We denote by $ \Yy_n $ the set of Young diagrams of size $n$ or, equivalently, the set of partitions of the integer $n$. We will set $ \lambda \vdash n $ for $ \lambda \in \Yy_n $. If $ \lambda \vdash n $, we set $ \abs{\lambda} = n $, we define $ \ell(\lambda) $ for its length, $ n(\lambda) := \sum_{k \geq 1} (k - 1) \lambda_k $ and $ m_k(\lambda) := \sum_{j \geq 1} \Unens{ \lambda_j  = k} $ for the multiplicity of $k$ in $ \lambda $ (see \cite[ch. I-1]{MacDo}). The set of strict partitions, i.e. partitions $ \lambda $ such that $ \lambda_1 > \lambda_2 > \dots $ will be denoted by $ \Yy_n^{(s)} $. Last, we denote by $ h(\square) $ the hook length of $ \square  \in \lambda  $.

\subsubsection{\textbf{Symmetric functions}}

We will use the following symmetric functions that are defined in \cite[ch. I]{MacDo} : the homogeneous complete symmetric functions $ h_n $, the power functions $ p_n $ and the Schur functions $ s_\lambda $.

Throughout this paper, we will make a constant use of the $ \lambda $-ring (or plethystic) notation (see e.g. \cite{Ram1991, LascouxSym}). The plethysm of a symmetric function $ f $ in a ``virtual'' alphabet $ \Ae $ will be denoted between brackets $ [ \, ] $, namely $ f[\Ae] $. We also define
\begin{eqnarray*}
H\crochet{\Ae} := \sum_{n \geq 0} h_n\crochet{\Ae} = \prod_{a \in \Ae} \frac{1}{1 - a} \ \ \ \mbox{if} \ \Ae = \ensemble{ a }_{a \in \Ae}
\end{eqnarray*}

The product of the alphabets $ \Ae $ and $ \Be $ is defined by $ \Ae \Be := \ensemble{ab}_{a \in \Ae, b \in \Be} $. We will use in particular the alphabet 
\begin{eqnarray*}
\frac{1}{1 - q} = \ensemble{ q^k }_{k \geq 0}
\end{eqnarray*}

For instance, for $ \abs{q} < 1 $,
\begin{eqnarray*}
p_n\crochet{ \frac{1}{1 - q} } = \sum_{k \geq 0} q^{kn} = \frac{1}{1 - q^n}
\end{eqnarray*}

\subsection{Reminders on mod-Poisson convergence}

 We define the mod-Poisson convergence (in the Laplace-Fourier setting) by :  

\begin{shaded}
\begin{definition}\label{Def:ModPoissonCv} 
Let $ (X_n)_n $ be a sequence of positive random variables, let $ (\gamma_n)_n $ be a sequence of strictly positive real numbers and $ \Peb(\gamma_n) \sim  \Poisson(\gamma_n) $. The sequence $ (X_n)_n $ is said to converge in the mod-Poisson sense at speed $ (\gamma_n)_n $ if there exists a continuous function $ \Phi : \Cc \to \Cc $ satisfying $ \Phi(1) = 1 $ and $ \Phi(\overline{z}) = \overline{\Phi(z) } $ such that the following convergence holds locally uniformly in $ z \in \Cc $
\begin{align*}
\frac{ \Esp{ z^{  X_n} } }{ \Esp{ z^{   \Peb(\gamma_n) } } } \tendvers{n }{ + \infty } \Phi\prth{ z }
\end{align*}

When such a convergence holds, we write it as
\begin{align*}
(X_n, \gamma_n ) \cvmodp{n }{ + \infty } \Phi
\end{align*}
\end{definition}	
\end{shaded}


\begin{remark} The limiting function $ \Phi $ is not unique. It is defined up to multiplication by an exponential due to the choice of speeds $ \gamma_n $ (see \cite{JacodAl}). 
\end{remark}

\medskip
One can also use some more restricted definition of mod-Poisson convergence, by only allowing $ z $ to belong to a subset of $ \Cc $, for instance the unit circle $ \Uu $ (the original definition) or taking $ z = e^{-\lambda} $, with $ \lambda > 0 $, i.e. $ z = x \in (0, 1] $ (Laplace setting).

\medskip\medskip
\section{The symmetric group}\label{Sec:Sn}

\subsection{Reminders on random uniform permutations}\label{SubSec:Sn:Reminders}

\subsubsection{\textbf{Cycle structure}}\label{SubSubSec:Sn:Reminders:Cycles}

The most classical example of mod-Poisson convergence is given by the total number of cycles of a random uniform permutation $ \sigmab_{\! n} \sim \Us(\Sg_n) $. Let $ C(\sigmab_{\! n}) $ denote this random variable. The classical chinese restaurant theorem (see e.g. \cite{ArratiaBarbourTavare}) implies that 
\begin{align*}
C(\sigmab_{\! n}) \eqlaw \sum_{k = 1}^n B_k, \qquad B_k \sim \Ber_{\{0, 1\}}\prth{\tfrac{1}{k}} \quad\mbox{(independent)}
\end{align*}

In particular, the Poisson approximation is valid, and one can write with $ \Peb(H_n) \sim \Poisson(H_n) $
\begin{align*}
\frac{ \Esp{ z^{ C(\sigmab_{\! n}) }  } }{\Esp{ z^{ \Peb(H_n) } } } = \frac{ \prod_{k = 1 }^n \prth{ 1 + \frac{z-1}{k} } }{ e^{ (z-1) \sum_{k = 1}^n \frac{1}{k} } }  = \prod_{k = 1 }^n \prth{ 1 + \frac{z-1}{k} } e^{ - \frac{z-1}{k} } \tendvers{n}{ + \infty } \prod_{k \geq 1} \prth{ 1 + \frac{z-1}{k} } e^{ - \frac{z-1}{k} } 
\end{align*}
the last limit being locally uniform in $ z \in \Cc $, for $  \prth{ 1 + \frac{z-1}{k} } e^{ - \frac{z-1}{k} } = \exp\prth{ - \frac{ (x-1)^2 }{2 k^2}  + o\prth{ \frac{1}{k^2} } } $ and $ \sum_k \frac{1}{k^2} < \infty $.

Taking instead $ \Peb(H_n - \gamma) $ or equivalently $ \Peb(\log n) $, we get 
\begin{align}\label{Thm:ModPoissonCvTotalNumberCycles}
\frac{ \Esp{ z^{ C(\sigmab_{\! n}) }  } }{\Esp{ z^{ \Peb(\log n) } } }  
\tendvers{n}{ + \infty }  
e^{\gamma (z-1) } \prod_{k \geq 1} \prth{ 1 + \frac{z-1}{k} } e^{ - \frac{z-1}{k} } 
\end{align}

\subsubsection{\textbf{Cycle Index Polynomial and Newton's Binomial Formula}}\label{SubSubSec:Sn:Reminders:CIP}

Newton's original binomial formula was the expansion of $ (1 - t)^{-x} $ which reads 
\begin{align*}
\frac{1}{(1 - t)^x } = \sum_{n \geq 0 } \frac{x^{\uparrow n} }{n!} t^n
\end{align*}
%
%
%

With $ g_t \sim \Geom(t) $ independent of $ (\sigmab_{\! n})_n $, one obtains
\begin{align*}
\Esp{ z^{ C(\sigmab_{\! g_t}) } } & 
                = \sum_{n \geq 0 } \Prob{ g_t = n } \Esp{ z^{ C(\sigmab_{\! n}) } } 
                = (1 - t) \sum_{n \geq 0 } \frac{z^{\uparrow n} }{n!} t^n 
                = \frac{1 - t}{(1 - t)^z } = e^{ (z - 1) \log\left( \frac{1}{1 - t} \right) } \\
                & = \Esp{ z^{ \Peb\left(\log\left( \frac{1}{1 - t} \right)\right) } }
\end{align*}
which implies, by injectivity of the Fourier-Laplace transform
\begin{align}\label{Thm:EqLawCsigmaGt}
C(\sigmab_{\! g_t}) \eqlaw    \Peb\prth{  \log\prth{\tfrac{1}{1 - t}} }
\end{align}

Newton's binomial formula thus has a strong probabilistic flavour since it expresses the fact that the total number of cycles of a geometrized uniform random permutation is Poissonian. This formula is a particular case of Polya's cycle index theorem that expresses the link between the symmetric power functions $ p_k(X) $ and the complete homogeneous functions $ h_n(X) $ (see e.g. \cite[ch. I-2 (2.10)]{MacDo}). Indeed, the Cycle Index Polynomial of $ \Sg_n $ is defined by
\begin{align*}
CIP(\Sg_n)(X) := \frac{ 1 }{n! } \sum_{ \sigma \in \Sg_n } \prod_{k \geq 1} p_k(X)^{c_k(\sigma)} 
               = \Esp{ \prod_{k = 1}^n p_k(X)^{ c_k(\sigmab_{\! n}) } } 
               = h_n(X)
\end{align*}
where $ c_k(\sigma) $ is the number of cycles of $ \sigma $ with length $ k $, and Polya's theorem reads
\begin{align*}
H(X) := \sum_{n \geq 0 } t^n h_n(X) = \exp\prth{ \sum_{k \geq 1 } \frac{p_k(X)}{k } t^k } =: e^{P(X)}
\end{align*}

Setting $ x_k := p_k(X) $ and multiplying by $ 1 - t =  \exp\prth{ - \sum_{k \geq 1} t^k/k  } $, one gets
\begin{align}\label{Thm:ProbabilisticPolyaCIP}
\Esp{\prod_{k \geq 1} x_k^{c_k(\sigmab_{\! g_t})}} 
                      = \sum_{n \geq 0 } \Prob{ g_t = k } \Esp{ \prod_{k = 1}^n x_k^{ c_k(\sigmab_{\! n}) } } 
                      = e^{ \sum_{k \geq 1 } \frac{t^k}{k } (x_k - 1)  } 
                      = \prod_{k \geq 1} \Esp{ x_k^{ \Peb(t^k/k) } }
\end{align}
which is equivalent by Fourier injectivity to
\begin{align}\label{Thm:EqLawCycleStructPoisson}
\prth{ c_k(\sigmab_{\! g_t}) \vphantom{\big)} }_{k \geq 1} \eqlaw \prth{ \Peb\prth{ \frac{t^k}{k} } }_{k \geq 1}
\end{align}
where the random variables $ \Peb(t^k/k) \sim \Poisson(t^k/ k) $ are independent. 


\begin{remark}
Using \eqref{Thm:EqLawCycleStructPoisson} and the fact that $ C(\sigma) = \sum_{k \geq 1} c_k(\sigma) $, it is not hard to deduce \eqref{Thm:EqLawCsigmaGt} using the independence and the infinite divisibility of the $ \Peb(t^k/k) $. Indeed, 
\begin{align*}
C(\sigma_{g_t}) = \sum_{k \geq 1} c_k(\sigma_{g_t}) 
                \eqlaw \sum_{k \geq 1} \Peb\prth{ \frac{t^k}{k} } 
                \eqlaw \Peb\prth{ \sum_{k \geq 1} \frac{t^k}{k} } 
                = \Peb\prth{ \log( (1 - t)\inv ) }
\end{align*}

In the same vein, using $ n = \sum_{k \geq 1 } k c_k(\sigmab_{\! n}) $, one gets the equality in law
\begin{align}\label{Thm:EqLawDimensionGeom}
\sum_{k \geq 1} k \Peb\prth{ \frac{t^k}{k} } 
                 \eqlaw \sum_{k \geq 1} k C_k(\sigma_{g_t}) 
                 = g_t \sim \Geom(t)
\end{align}

The fact that the geometric distribution of parameter $ t $ can be written as an infinite linear combination of Poisson random variables is the probabilistic translation of the expansion $ \log\prth{1 - t} = - \sum_{k \geq 1} t^k / k $ (taking the Laplace transform of this equality in law).
\end{remark}

\subsection{Mod-Poisson convergence of $ C(\sigmab_{\! n}) $ using the geometric randomisation}\label{SubSec:Sn:ModPo}

We now prove \eqref{Thm:ModPoissonCvTotalNumberCycles} using $ g_t $. This proof is the prototype of all future proofs that use the randomisation paradigm at stake in this article to explain the appearance of the $ \Gamma $ factor.

\begin{shaded}
\begin{theorem}\label{theorem:FacteurGammaCycles} Suppose that $ (C(\sigmab_{\! n}))_n $ converges in the mod-Poisson sense at speed $ \ln(n) + \gamma $ to a certain function $ \Phi $ with the following speed of convergence:
\begin{align*}
\frac{ \Esp{ x^{ C(\sigmab_{\! n}) } } }{  \Esp{ x^{\Peb(\log n + \gamma)} } } 
                 = \Phi(x) \prth{ 1 + O_x\prth{ \frac{1}{n} } } 
\end{align*}

Then, necessarily, $ \Phi(x) = \frac{ 1 }{\Gamma(x) } $.
\end{theorem}
\end{shaded}


\begin{proof}
Define $ \Phi_\gamma(x) := \Phi(x)e^{\gamma(x - 1)} $. By hypothesis, we have 
\begin{align*}
\Esp{ x^{ C(\sigmab_{\! n}) } } 
                 = \Phi(x)  \Esp{ x^{\Peb(\log n + \gamma)} }(1 + O_x(n\inv) ) 
                 = \Phi_\gamma(x) n^{x - 1} (1 + O_x(n\inv) ) 
\end{align*}

Let us \textit{randomise} this last equality, writing it for $ n = g_t $ with $ t \to 1 $ (i.e. $ n \to +\infty $). We have
\begin{align*}
\Esp{ x^{ C(\sigmab_{\! g_t}) } } 
                   = \Phi_\gamma(x) \Esp{ g_t^{x - 1}  + O_x\prth{ g_t^{x - 2} } } 
                   = \Phi_\gamma(x) \crochet{ \Esp{ g_t^{x - 1} } + O_x\prth{ \Esp{ g_t^{x - 2} } }  }
\end{align*}
since the constant in the $ O_x $ is deterministic.

Using $ g_t = \pe{ \ee / \log(t\inv) } 
            = \ee / \log(t\inv) - \pf{ \ee / \log(t\inv) } 
            = \ee / \log(t\inv) + O_{t, \ee}(1)  
            = \ee / \log(t\inv) + O(1) $ since $ 0 \leq \pf{ \ee / \log(t\inv) } < 1 $ a.s., we get 
\begin{align*}
\Esp{ g_t^{x - 1} } & =  \Esp{ \prth{\ee / \log(t\inv) + O(1) }^{x - 1} }  \\
                    & =  (1/\log(t\inv))^{x - 1}  \Esp{ \prth{\ee   + O(\log(t\inv)) }^{x - 1} }   \\
                    & =  e^{(x - 1) \log( 1/\log(t\inv) )}  \Esp{ \ee^{x - 1} } \prth{ 1 + O_x\prth{ (\log(t\inv))^{x-1}} } 
\end{align*}
%
%
%
%
%
%
and in particular, $ O_x\prth{ \Esp{ g_t^{x - 2} } } = O_x\prth{ (1/\log(t\inv))^{x - 2} } $. We will now set $ t := 1 - \varepsilon $ with $ \varepsilon \to 0 $; we then have $ \log(t\inv) = \log(\frac{1}{1 - \varepsilon}) = \varepsilon + O(\varepsilon^2) $. From this, we deduce that, locally uniformly in $ x \in \Cc $
\begin{align*}
\Esp{ x^{ C(\sigmab_{\! g_{1 - \varepsilon}}) } } 
                & = \Phi_\gamma(x) 
                     e^{(x - 1) \log( \varepsilon\inv )(1 + o(1)) }  
                     \Esp{ \ee^{x - 1} } \prth{ 1 + O\prth{ \varepsilon } } \\
                & = \Phi_\gamma(x) 
                    \Gamma(x) 
                    e^{(x - 1) \log(\varepsilon\inv) } \prth{ 1 + O\prth{ \varepsilon } } \quad \mbox{with \eqref{Eq:TransfoMellinExp}}
\end{align*}

We thus get on the one hand
\begin{align*}
\frac{\Esp{ x^{ C(\sigmab_{\! g_{1 - \varepsilon} }) } } }{ e^{(x - 1) \log(\varepsilon\inv) } } 
               =  \Phi_\gamma(x) \Gamma(x) \prth{ 1  + O\prth{  \varepsilon } }  
               \tendvers{t}{1} \Phi_\gamma(x) \Gamma(x) 
\end{align*}

On the other hand, let $ v_t := \log\prth{ \frac{1}{1 - t } } = \ln(\varepsilon\inv) $. By \eqref{Thm:EqLawCsigmaGt}, we have $ C(\sigmab_{\!g_t}) \eqlaw \Peb(v_t) $ and in particular, $ \Esp{ x^{ C(\sigmab_{\! g_t}) } } = e^{(x-1) v_t } $, yielding 
\begin{align*}
\frac{ \Esp{ x^{ C(\sigmab_{\! g_t}) } } }{ e^{ (x-1) \log(1/\log(t\inv) ) } }  
              = e^{(x-1) (v_t + \log\log (t\inv) ) } 
              =  e^{(x-1) \log\left( \frac{ -\log(1 - \varepsilon)}{ \varepsilon}  \right) }
              \tendvers{t}{1^-} 1 
\end{align*}
%
%
%
%
%
%
which implies that $ \Phi_\gamma(x) \Gamma(x) = 1 $. 
\end{proof}


\begin{shaded}
\begin{proposition}\label{proposition:ConvergenceModPoissonCycles} 
The convergence \eqref{Thm:ModPoissonCvTotalNumberCycles} is satisfied.
\end{proposition}
\end{shaded}


\begin{proof}
We need to show that mod-Poisson convergence exists at speed $ \log n $ with a speed of convergence towards the limiting function equal to $ O_x(1/n) $. One has 
\begin{align*}
\frac{\Esp{ x^{ C(\sigmab_{\! n}) } } }{ e^{(x - 1) \log( n )} }  & = e^{(x-1)(H_n - \log n) } \exp\prth{ \sum_{k = 1}^n \log\prth{ 1 + \frac{x-1}{k} } - \frac{x-1}{k} } \\
                       & = \exp\prth{ \prth{ \gamma + O\prth{ \frac{1}{n} } }(x-1) + \sum_{k \geq 1}  \log\prth{ 1 + \frac{x-1}{k} } - \frac{x-1}{k} +  (x-1)^2 O\prth{\frac{1}{ n }  } }
\end{align*}
as $  \log\prth{ 1 + \frac{x-1}{k} } - \frac{x-1}{k}  = O\prth{  \frac{(x-1)^2 }{ k^2 } } $ for $ x - 1 \in \Cc\setminus\Rr_-  $  and $ \sum_{ k \geq n+1 } k^{-2 } = O(1/n) $. We thus have
\begin{align*}
\frac{\Esp{ x^{ C(\sigmab_{\! n}) } } }{ e^{(x - 1) \log( n )} }  & = \exp\prth{ \gamma (x-1) + \sum_{k \geq 1}  \log\prth{ 1 + \frac{x-1}{k} } - \frac{x-1}{k} +    O_x\prth{\frac{1}{ n }  } } \\
                       & = \exp\prth{ \gamma (x-1) + \sum_{k \geq 1}  \log\prth{ 1 + \frac{x-1}{k} } - \frac{x-1}{k} }\prth{ 1 +    O_x\prth{\frac{1}{ n }  } }
\end{align*}
hence the result.
\end{proof}


\begin{remark}\label{Rk:BigOvsLittleO}
Replacing $ O_x(n\inv) $ by $ o_x(1) $ in theorem~\ref{theorem:FacteurGammaCycles} also gives the result. The only reason we supposed such a strong speed of convergence comes from proposition~\ref{proposition:ConvergenceModPoissonCycles} and from the desire to explain which kind of speed one gets with the randomisation. 
\end{remark}


\begin{remark}
Another way to conclude is given by Shepp and Lloyd \cite[(4)]{SheppLloyd1966} who use a ``Binet form of the Stirling approximation'' for $ z \in \Uu $, citing \cite[p. 249]{WhittakerWatson}.
\end{remark}


\begin{remark}\label{Rk:EulerWithMod}
We have thus proven that 
\begin{align}\label{Eq:EulerFormulaGamma}
\frac{1}{\Gamma(x) } = e^{\gamma(x - 1) } \prod_{k \geq 1} \prth{ 1 + \frac{x-1}{k} } e^{ -\frac{x-1}{k} }
\end{align}
since the RHS is the limiting function at speed $ \log n + \gamma $ given in \eqref{Thm:ModPoissonCvTotalNumberCycles} and the LHS is the one we obtain with theorem \ref{theorem:FacteurGammaCycles} at speed $ \log n + \gamma $. This formula was first discovered by Euler, see e.g. \cite{WhittakerWatson}. The last computations hence give a probabilistic proof of this identity.

Another probabilistic interpretation of this identity goes as follows: write it
\begin{align*}
\Esp{ \ee^{x-1} } = \Gamma(x)  = e^{- \gamma(x - 1) } \prod_{k \geq 1} \frac{1}{ 1 + \frac{x-1}{k} } e^{  \frac{x-1}{k} } = e^{- \gamma(x - 1) } \prod_{k \geq 1} \Esp{ e^{  \frac{x-1}{k} \ee_k } } e^{  \frac{x-1}{k} }
\end{align*}
where $ (\ee_k)_k $ is a sequence of i.i.d. random variables with $ \Exp(1) $ distribution. The Euler identity is thus equivalent to the equality in law (whose RHS clearly converges in $ L^2 $)
\begin{align*}
\log \ee \eqlaw - \gamma + \sum_{k \geq 1} \frac{ \ee_k - 1 }{k }
\end{align*}
\end{remark}

\medskip\medskip
\section{The Sathe-Selberg theorems and the delta-Zeta distribution}\label{Sec:Integers}

\subsection{Analogies between permutations and integers}\label{SubSec:Integers:Analogies}

We wish to extend the previous setting to number fields in order to prove \eqref{Thm:SatheSelberg} with a randomisation. To achieve this goal, we need an equivalent of the geometric randomisation at stake in Polya's CIP theorem \eqref{Thm:ProbabilisticPolyaCIP} or Newton's formula \eqref{Thm:EqLawCsigmaGt}.

There are several analogies between permutations and integers ; some are listed in \cite{Grandville} or in \cite[p. 22]{ArratiaBarbourTavare}. The equivalent of the cycle decomposition of $ \sigma \in \Sg_n $ is the prime decomposition \eqref{Thm:PrimeDecomposition} of $ k \in \Nn^* $ ; in this analogy, the equivalent of the number of $m$-cycles $ c_m(\sigma) $ is the $p$-adic valuation $ v_p(k) $. If we assume that the analogue of the uniform measure $ \Us(\Sg_n) $ is $ \Us(\intcrochet{1, n}) $, we thus want the equivalent $ d_t $ of the randomisation $ g_t \sim \Geom(t) $ such that the valuations $ (v_p(U_{d_t}))_p $ are independent (with $ U_n \sim \Us(\intcrochet{1, n}) $). This distribution is simpler than the ``geometrised uniform distribution on $ \bigcup_{n \geq 0} \Sg_n $'' as it is still a distribution on $ \Nn^* $.

The first question to ask is therefore : is there a \textit{natural} random variable with values in $ \Nn^* $ whose $p$-adic valuations are independent ?

\medskip
\subsection{Reminders on the Zeta distribution}\label{SubSec:Integers:Zeta}

Recall that the Riemann $ \zeta $ function is defined for $ \alpha \in \ensemble{ \Re > 1 } $ by its Dirichlet series or its Euler product
\begin{align*}
\zeta(\alpha) := \sum_{n \geq 1 } \frac{1 }{n^\alpha } = \prod_{p \in \Pe} \frac{1 }{ 1 - p^{-\alpha } }
\end{align*}

\begin{shaded}
\begin{definition}[Zeta distribution] Let $ \alpha > 1 $. The Zeta distribution of parameter $ \alpha $ is denoted by $ \zetab(\alpha) $ and defined, with $ Z_\alpha \sim \zetab(\alpha) $, by
\begin{align}\label{Def:ZetaDistribution}
\forall k \geq 1, \ \ \ \  \Prob{ Z_\alpha = k } = \frac{1}{ \zeta(\alpha) k^\alpha }
\end{align}

This is clearly a probability on $ \Nn^* $ by definition of the Riemann $ \zeta $ function.
\end{definition}
\end{shaded}

The main property\footnote{
As the geometric distribution $ \Geom(t) $ is infinitely divisible, of L\'evy measure $ \sum_{k \geq 1} \frac{t^k}{k} \delta_k $, $ \zetab(\alpha) $ is the exponential of an infinitely divisible distribution.
} of this distribution is indeed to have independent $p$-adic valuations :

\begin{lemma}[Khintchin, \cite{Khintchin}] 
\begin{shaded}
Let $ \Geb\prth{ p^{ -\alpha } } \sim \Geom\prth{ p^{ -\alpha } } $ be independent random variables and let $ Z_\alpha \sim \zetab(\alpha) $ for $ \alpha > 1 $. Then,
\begin{align}\label{Thm:KintchineEqLaw}
Z_\alpha \eqlaw \prod_{p \in \Pe }  p^{ \Geb(p^{ -\alpha } ) }
\end{align}

In particular, $ (v_p(Z_\alpha))_{p \in \Pe} $ are independent random variables with a geometric distribution:
\begin{align}\label{Thm:pValuationsIndependance}
v_p(Z_\alpha) \eqlaw \Geb\prth{ p^{ -\alpha } }  \sim \Geom\prth{ p^{ -\alpha } }  
\end{align}
\end{shaded}
\end{lemma}

For the reader's convenience, we recall the proof of \eqref{Thm:KintchineEqLaw}.  

\begin{proof}
Let $ \xb :=(x_p)_{p \in \Pe} $ with $ x_p < p^\alpha $ and set $ \vb := (v_p)_{p \in \Pe} $. Then,
\begin{align*}
\Esp{ \xb^{ \vb(Z_\alpha) } } = \frac{1}{\zeta(\alpha) } \sum_{n \geq 1 } \frac{1}{n^\alpha } \prod_{p \in \Pe } x_p^{ v_p(n) } = \prod_{p \in \Pe } \prth{ 1 - \frac{1}{p^\alpha } } \sum_{n \geq 1 }  \prod_{p \in \Pe } \prth{ \frac{ x_p}{p^\alpha } }^{ v_p(n) }
\end{align*}

Let $ (a_p)_p $ with $ \abs{a_p} < 1 $. Supposing that all considered series and products are convergent, we have the following Euler's identity
\begin{align}\label{Thm:FormuleEulerProduitSommePremiers}
\sum_{n \geq 1 }  \prod_{p \in \Pe } a_p^{ v_p(n) } = \prod_{p \in \Pe } \frac{1}{1 - a_p}
\end{align}

Indeed, 
\begin{align*}
\prod_{p \in \Pe } \frac{1}{1 - a_p} 
                  & = \prod_{p \in \Pe } \sum_{k_p \geq 0 } a_p^{k_p } 
                    = \sum_{ (k_p)_{p \in \Pe } \in \Nn^\Pe } \prod_{p \in \Pe } a_p^{k_p } 
                    = \sum_{n \geq 1 } \sum_{ (k_p)_{p \in \Pe } } \Unens{ (k_p)_p = \vb(n) } \prod_{p \in \Pe } a_p^{ v_p(n) } \\
                 & = \sum_{n \geq 1 }  \prod_{p \in \Pe } a_p^{ v_p(n) }
\end{align*}

We deduce from \eqref{Thm:FormuleEulerProduitSommePremiers} with $ a_p := x_p/p^\alpha < 1 $ that
\begin{align*}
\Esp{ \xb^{ \vb(Z_\alpha) } }  
             =  \prod_{p \in \Pe } \frac{1 - \frac{1}{p^\alpha }}{ 1 - \frac{ x_p}{p^\alpha } } 
             = \prod_{p \in \Pe } \Esp{ x_p^{ \Geb( p^{-\alpha} ) } }
\end{align*}
which gives the desired result.
\end{proof}

\begin{remark}
Khintchin's equality in law \eqref{Thm:KintchineEqLaw} is a particular case of a most general type of distribution, the multiplicative distributions. These are precisely defined by an equality in law of the type $ \prod_{p \in \Pe} p^{X_p} $ where $ (X_p)_p $ is a sequence of independent random variables. The $p$-adic valuations of these random variables are thus, by construction, independent, but $ \Prob{ \prod_{p \in \Pe} p^{X_p} = n } = \prod_{p \in \Pe } \Prob{ X_p = v_p(n) }  $ is not as easily expressible as in the case of $ \zetab(\alpha) $. 
\end{remark}

We conclude these reminders with the fluctuations of $ \zetab(\alpha) $ when $ \alpha \to 1^+ $.

\begin{shaded}
\begin{lemma} Let $ Z_\alpha \sim \zetab(\alpha) $ and $ \ee \sim \Exp(1) $. Then, 
\begin{align}\label{Thm:FluctuationsLoiZeta}
( \alpha - 1 ) \log Z_\alpha \cvlaw{\alpha }{ 1^+ } \ee
\end{align}
\end{lemma}
\end{shaded}


\begin{proof}
Comparing the sum to an integral, we have \cite[ch. I.1.6, proof thm. 11 p. 17]{Tenenbaum}
\begin{align}\label{Thm:ComportementZetaEn1}
\zeta(1 + \varepsilon) = \frac{1}{\varepsilon} + O(1) \ \ \ \mbox{when } \varepsilon \to 0^+
\end{align}

As a result,
\begin{align*}
\Prob{ Z_\alpha \leq r }  & = \frac{1}{\zeta(\alpha) } \sum_{1 \leq k \leq r} \frac{1}{k^\alpha }  \\
                & = \frac{1}{ \zeta(\alpha) } \prth{ \int_1^r \frac{dx}{ x^\alpha } + O\prth{ \frac{1}{ r^{\alpha - 1} } } } \\
                & = \frac{1}{\zeta(\alpha) } \prth{  \frac{1 - r^{1 - \alpha} }{\alpha - 1 } + O\prth{ \frac{1}{ r^{\alpha - 1} } } } \\  
                & =  \frac{1 - r^{1 - \alpha} }{ 1 + o(1) } + O\prth{ \frac{\alpha - 1}{   r^{\alpha - 1} } } \quad\mbox{with \eqref{Thm:ComportementZetaEn1}.} 
\end{align*}

Taking $ r = \pe{ \exp\prth{ t / (\alpha - 1) } } = \exp\prth{ t / (\alpha - 1) } + O(1) $, we thus get
\begin{align*}
\Prob{ (\alpha - 1) \log Z_\alpha \leq  t  }  =  \frac{1 - e^{ - t } }{ 1 + o(1) } + O\prth{ (\alpha - 1) e^{ -t  } } \tendvers{ \alpha }{ 1 } 1 - e^{-t } = \Prob{ \ee \leq t }
\end{align*}
which gives the desired result.
\end{proof}


\medskip
\subsection{The delta-Zeta distribution}\label{SubSec:Integers:DeltaZeta}

\subsubsection{\textbf{Definition and properties}}\label{SubSubSec:Integers:DeltaZeta:Def}

We now introduce a new distribution on $ \Nn^* $ such that the ``arithmetic cycle index theorem'' (in probabilistic form)
\begin{align}\label{Thm:ProbabilisticArithmeticCIP}
\sum_{n \geq 1} \Prob{ d_t = n } \Esp{ \xb^{\vb(U_n) } } = \prod_{p \in \Pe} \Esp{ x_p^{ \Geb( p^{-t } ) } }
\end{align}
becomes the direct analogue of Polya's cycle index theorem in probabilistic form \eqref{Thm:ProbabilisticPolyaCIP}
\begin{align*}
\sum_{n \geq 0} \Prob{ g_t = n } \Esp{ \xb^{\cb(\sigmab_{\! n}) } } = \prod_{k \in \Nn^*} \Esp{ x_k^{ \Peb( t^k/k ) } }
\end{align*}


\begin{shaded}
\begin{definition}[Delta-Zeta distribution] The delta-Zeta distribution of parameter $ \alpha > 1 $ is denoted by $ \deltab\zetab(\alpha) $ and defined, with $ d_\alpha \sim \deltab\zetab(\alpha) $, by
\begin{align}\label{Def:DeltaZetaDistribution}
\forall k \geq 1, \ \ \ \  \Prob{ d_\alpha = k } = \frac{k}{ \zeta(\alpha)} \prth{ \frac{1}{ k^\alpha } - \frac{1}{ (k + 1)^\alpha } }
\end{align}
\end{definition}
\end{shaded}

It is clear that $ \sum_{k \geq 1} \Prob{ d_\alpha = k } = 1 $ by an Abel summation.

\medskip

Let $ U \sim \Us([0, 1]) $ and $ R_\alpha := \pe{ U^{1/\alpha } } $. Then, 
\begin{align*}
\Prob{R_\alpha = k } = \frac{1}{ k^\alpha } - \frac{1}{ (k + 1)^\alpha } \ \  \ \ \mbox{and } \ \  \ \ \Esp{R_\alpha } = \zeta(\alpha)
\end{align*}

Thus, $ d_\alpha \sim \deltab\zetab(\alpha) $ is the size-biased transform of $ R_\alpha $ :
\begin{align}\label{Eq:DeltaZeta=SizeBias}
\Prob{ d_\alpha = k } = \frac{k}{ \Esp{R_\alpha } } \Prob{ R_\alpha = k }  
\end{align}

We list some additional properties of $ \deltab\zetab(\alpha) $ in Annex~\ref{Sec:DeltaZeta}.

\medskip
\subsubsection{\textbf{The arithmetic CIP through Delta-Zeta randomisation}}\label{SubSubSec:Integers:DeltaZeta:ArithmeticCIP}

We now prove \eqref{Thm:ProbabilisticArithmeticCIP}.


\begin{shaded}
\begin{theorem}[Delta-Zeta randomisation] 
Let $ n \geq 1 $ and $ \alpha > 1 $ ; let $ U_n \sim \Us\prth{ \intcrochet{1, n} } $, $ d_\alpha \sim \deltab\zetab(\alpha) $ and $ Z_\alpha \sim \zetab(\alpha) $, where $ (U_n)_n $ and $ d_\alpha $ are independent. We denote by $ U_{d_\alpha } $ the independent randomisation given by $ \Prob{ U_{d_\alpha } = k } = \sum_{n \geq 1} \Prob{ U_n = k } \Prob{ d_\alpha = n } $. Then, 
\begin{enumerate}
\item the following equality in law is satisfied
\begin{align}\label{Thm:RandomisationDeltaZeta}
U_{d_\alpha } \eqlaw Z_\alpha 
\end{align}

\item and this equality in law defines uniquely the law of $ ( d_\alpha, Z_\alpha ) $ in the sense that if an independent randomisation $ X $ of a uniformly distributed random variable gives a random variable $ Y $ with independent geometric $p$-adic valuations, then, there exists a unique $ \alpha > 1 $ such that $ X \sim \deltab\zetab(\alpha) $ and $ Y \sim \zetab(\alpha) $.
\end{enumerate}
\end{theorem}
\end{shaded}


\medskip
\begin{proof}

$ $

\noindent $ (1) $ Let $ s \in \Rr $. Then, 
\begin{align*}
\Esp{ U_{d_\alpha}^{is} } & = \sum_{k \geq 1 } \Prob{ d_\alpha = k } \Esp{  U_k^{is} } \\
						  & =  \frac{1}{\zeta(\alpha) } \sum_{k \geq 1} k \prth{ \frac{1}{k^\alpha} - \frac{1}{ (k + 1)^\alpha } } \frac{1}{k} \sum_{\ell = 1}^k \ell^{is} 
						    = \frac{1}{\zeta(\alpha) } \sum_{k \geq \ell \geq 1} \prth{ \frac{1}{k^\alpha} - \frac{1}{ (k + 1)^\alpha } }\ell^{is} \\
						  & = \frac{1}{\zeta(\alpha) } \sum_{\ell \geq 1} \frac{ 1 }{\ell^\alpha} \ell^{is} \\
						  & = \Esp{ Z_\alpha^{is} } 
\end{align*}

We conclude by injectivity of the Fourier-Mellin transform.

\medskip

\noindent $ (2) $ Let $ (X, Y) $ satisfying the subordination equation $ U_X \eqlaw Y $ with $ Y \eqlaw \prod_{p \in \Pe} p^{ W_p } $ and $ (W_p)_p $ independent with a geometric distribution : $ W_p \sim \Geom(w_p) $. We define $ f : \Nn^* \to \Rr_+ $ by
\begin{align*}
\Prob{ Y = k } := \frac{f(k)}{ \sum_{ \ell \geq 1} f(\ell) }
\end{align*}

A fundamental property of $f$ is that it is a \textit{completely multiplicative function}, namely, $ f\prth{ \prod_{p \in \Pe} p^{v_p(n) } } = \prod_{p \in \Pe} f(p)^{v_p(n) } $. Indeed, for all $ s \in \Rr $,
\begin{align*}
\Esp{ Y^{is} } & = \frac{1}{ \sum_{ \ell \geq 1} f(\ell) } \sum_{k \geq 1} f(k) k^{is}  = \prod_{p \in \Pe } \Esp{ p^{i s W_p } } = \prod_{p \in \Pe} \frac{1 - w_p}{1 - w_p p^{is} } \\
               & =  \prod_{p \in \Pe} \prth{1 - w_p} \sum_{k \geq 1} k^{is} \prod_{p \in \Pe } w_p^{ v_p(k) }
\end{align*}
where all considered quantities are supposed to converge. By identification of the coefficient of $ k^{is} $ (by linear independence of the vector basis $ (s \mapsto e^{ i s \log k })_{k \geq 1} $), we deduce for $ k = 1 $ that $ \sum_{ \ell \geq 1} f(\ell) = \prod_{p \in \Pe} \prth{1 - w_p}\inv $ and for all $ k \geq 1 $
\begin{align*}
f(k) = \prod_{p \in \Pe } w_p^{ v_p(k) } = \prod_{p \in \Pe } f(p)^{ v_p(k) }
\end{align*}

Now, if $ U_X \eqlaw Y $, then
\begin{align*}
\Esp{ U_{X}^{is} } & =  \sum_{k \geq 1 } \Prob{X = k } \frac{1}{k} \sum_{\ell = 1}^k \ell^{is} =  \sum_{k \geq 1 } k^{is}\prth{ \sum_{\ell \geq k} \frac{1}{ \ell } \Prob{ X = \ell }  } \\
\Esp{ Y^{is} }    & =  \sum_{k \geq 1 } k^{is} \frac{f(k)}{\sum_{\ell \geq 1} f(\ell) } 
\end{align*}
which implies that for all $ k \geq 1 $, 
\begin{align*}
\sum_{\ell \geq k} \frac{1}{ \ell } \Prob{ X = \ell }  = \frac{f(k)}{\sum_{\ell \geq 1} f(\ell) }
\end{align*}

Solving this linear equation in $ \Prob{ X = k } $ gives
\begin{align*}
\Prob{ X = k }  = \frac{k}{\sum_{\ell \geq 1} f(\ell) } \prth{ f(k) - f(k + 1) }
\end{align*}

It is clear by Abel summation that 
\begin{align*}
\sum_{k \geq 1} \frac{k}{\sum_{\ell \geq 1} f(\ell) } \prth{ f(k) - f(k + 1) } = 1
\end{align*}
so one only needs to have $ f(k) - f(k + 1) \geq 0 $ in order to have a probability distribution. But a classical theorem on extremal orders of arithmetic functions (see e.g. \cite[pp. 82]{Tenenbaum}) characterizes the only multiplicative decreasing functions on $ \Nn^* $ as being the power functions, i.e. there exists a unique $ \alpha \in \Rr_+ $ such that
\begin{align*}
f(k) = k^{-\alpha}
\end{align*}

Finally, the convergence of $ \sum_k f(k) $ implies that $ \alpha > 1 $, concluding the proof.
\end{proof}

\medskip
\subsubsection{\textbf{Fluctuations of $ \deltab\zetab(\alpha) $ when $ \alpha \to 1 $}}\label{SubSubSec:Integers:DeltaZeta:Fluctuations}

\begin{shaded}
\begin{lemma} Let $ d_\alpha \sim \deltab\zetab(\alpha) $ and $ \ee \sim \Exp(1) $. Then, 
\begin{align}\label{Thm:FluctuationsLoiDeltaZeta}
( \alpha - 1 ) \log d_\alpha \cvlaw{\alpha }{ 1^+ } \ee
\end{align}
\end{lemma}
\end{shaded}


\begin{proof}
We have, using an Abel summation
\begin{align*}
\Prob{ d_\alpha \leq r }  & =  \frac{1}{\zeta(\alpha)} \sum_{k = 1}^r k \prth{ \frac{1}{k^\alpha} - \frac{1}{(k + 1)^\alpha } } = \frac{1}{\zeta(\alpha)} \sum_{k = 1}^r  \frac{1}{k^\alpha} -  \frac{1}{\zeta(\alpha) (r+1)^\alpha } \\
                          & = \Prob{ Z_\alpha \leq r }  + O\prth{ \frac{\alpha - 1 }{ (r+1)^\alpha } }
\end{align*}

Setting $ r = \pe{ e^{ t/ (\alpha - 1) } }  $ and using \eqref{Thm:FluctuationsLoiZeta} gives then the result. 
\end{proof}

\medskip
\subsection{Splitting of $ \Phi_\omega $ using the delta-Zeta randomisation}\label{SubSec:Integers:ModPoissonOmegaIntegers}

We now give the prototype of proof of universality of the Gamma factor using the (fluctuations of the) auxiliary randomisation, here $ d_\alpha \sim \deltab\zetab(\alpha) $. The ingredients of the proof are always the same, as a result, we will mainly focus on $ \omega(U_n) $ in this part, leaving the case of $ \Omega(U_n) $ to the reader.


\begin{shaded}
\begin{theorem}\label{theorem:FacteurGammaSatheSelberg} 
Suppose that $ (\omega(U_n))_n $ (resp. $ (\Omega(U_n))_n $) converges in the mod-Poisson sense at speed $ \log \log n + \kappa $ for a certain constant $ \kappa $ to a certain function $ \Phi_\omega $ (resp. $ \Phi_\Omega $) with a speed of convergence given by
\begin{align}\label{Thm:SatheSelberg:bis}
\begin{aligned}
\frac{ \Esp{ x^{ \omega(U_n) } } }{  \Esp{ x^{\Peb(\log \log n + \kappa )} } } & = \Phi_\omega(x) \prth{ 1 + O_x\prth{ \frac{1}{\log n} } } \\
\frac{ \Esp{ x^{ \Omega(U_n) } } }{  \Esp{ x^{\Peb(\log \log n + \kappa )} } } & = \Phi_\Omega(x) \prth{ 1 + O_x\prth{ \frac{1}{\log n} } } 
\end{aligned}
\end{align}

Denote by $ \Phi_{ \widehat{\omega} }  $ (resp. $ \Phi_{ \widehat{\Omega} } $) the limiting function of the independent model $ \widehat{\omega}_n $ (resp. $\widehat{\Omega}_n $) defined in \eqref{Def:IndependentModel}.

Then, necessarily, there exists $ \widetilde{\kappa} \in \Rr $ such that $ \Phi_\omega(x) = \frac{ e^{\widetilde{\kappa} (x - 1)} }{ \Gamma(x) } \Phi_{ \widehat{\omega} }(x)  $ (resp. $ \Phi_\Omega(x) = \frac{ e^{\widetilde{\kappa} (x - 1)} }{\Gamma(x) } \Phi_{ \widehat{\Omega} }(x) $).
\end{theorem}
\end{shaded}


\begin{proof}
Let $ Z_\alpha \sim \zetab(\alpha) $. Since the $p$-adic valuations of $ Z_\alpha $ are independent and geometrically distributed, one has 
\begin{align*}
\omega(Z_\alpha) & = \sum_{p \in \Pe } \Unens{ v_p(Z_\alpha) \geq 1 } 
                   \eqlaw \sum_{ p \in \Pe } \Beb\prth{ p^{-\alpha} }, 
\qquad  \Beb\prth{ p^{-\alpha} } \sim \Ber_{\{0, 1\}} \prth{ p^{-\alpha} } \\
\Omega(Z_\alpha) & = \sum_{p \in \Pe }  v_p(Z_\alpha)  
                   \eqlaw \sum_{ p \in \Pe } \Geb\prth{ p^{-\alpha} }, 
\qquad \Geb\prth{ p^{-\alpha} } \sim \Geom\prth{ p^{-\alpha} }
\end{align*}
where all sequences are made of independent random variables.

We now focus on the case of $ \omega(U_n) $. By independence of the Bernoulli random variables, we have, with an absolutely convergent product in the RHS
\begin{align*}
\Esp{ x^{ \omega(Z_\alpha) } } = \prod_{p \in \Pe } \prth{ 1 + \frac{x-1}{p^\alpha} }
\end{align*}

For $ \alpha > 1 $, we define the \textit{prime $ \zeta $ function} by
\begin{align}\label{Def:PrimeZetaFunction}
\zeta_\Pe(\alpha) := \sum_{p \in \Pe } \frac{1}{p^\alpha}
\end{align}

Then, 
\begin{align*}
\Esp{ x^{ \omega(Z_\alpha) } } = e^{ (x-1) \zeta_\Pe(\alpha) } \prod_{p \in \Pe } \prth{ 1 + \frac{x-1}{p^\alpha} } e^{-\frac{x-1}{p^\alpha} } 
\end{align*}
and in particular, locally uniformly in $ x \in \Cc $, 
\begin{align*}
\frac{\Esp{ x^{ \omega(Z_\alpha) } }}{ \Esp{ x^{ \Peb(\zeta_\Pe(\alpha)) } } } =  \prod_{p \in \Pe } \prth{ 1 + \frac{x-1}{p^\alpha} } e^{-\frac{x-1}{p^\alpha} }  \tendvers{\alpha }{1 } \Phi_{ \widehat{\omega} } (\alpha) 
\end{align*}

By hypothesis and setting $ \Phi_{\omega, \kappa}(x) := \Phi_\omega(x) e^{\kappa (x - 1)} $, we have
\begin{align*}
\Esp{ x^{ \omega(U_n) } } 
                & = \Phi_\omega(x) \Esp{ x^{\Peb(\log \log n + \kappa)} } \prth{ 1 + O_x\prth{ \frac{1}{\log n} } } \\
                & = \Phi_{\omega, \kappa}(x)  (\log n)^{x - 1}  \prth{ 1 + O_x\prth{ \frac{1}{\log n} } }
\end{align*}

Let us randomise this equality by setting $ n = d_\alpha \sim \deltab\zetab(\alpha) $. On the one hand, we have 
\begin{align*}
\Esp{ x^{ \omega(U_{d_\alpha}) } } & = \Phi_{\omega, \kappa}(x)  \crochet{ \Esp{ (\log d_\alpha)^{x - 1} } + O_x\prth{ \Esp{ (\log d_\alpha)^{x - 2} }  } \vphantom{\Big)} } 
\end{align*}

Multiplying by $ (\alpha - 1)^{x-1} $ and using \eqref{Thm:FluctuationsLoiDeltaZeta} and \eqref{Eq:TransfoMellinExp}, we get
\begin{align*}
(\alpha - 1)^{x-1} \Esp{ x^{ \omega(U_{d_\alpha}) } } 
                 & = \Phi_{\omega, \kappa}(x)  \crochet{ \Esp{ \crochet{ (\alpha - 1) \log d_\alpha  }^{x - 1} } + \Esp{ \crochet{ (\alpha - 1) \log d_\alpha  }^{x - 2} } O_x\prth{ \alpha - 1 } } \\
                 & \tendvers{\alpha }{ 1 } \Phi_{\omega, \kappa}(x) \Esp{ \ee^{x-1} } 
                   = \Phi_{\omega, \kappa}(x) \Gamma(x) 
\end{align*}

On the other hand, using \eqref{Thm:RandomisationDeltaZeta}, we have 
\begin{align*}
(\alpha - 1)^{x-1} \Esp{ x^{ \omega(U_{d_\alpha}) } } 
                   & = \frac{ \Esp{ x^{ \omega(Z_\alpha) } } }{ \Esp{ x^{ \Peb( -\log (\alpha - 1)  ) } } } 
                   = \frac{ \Esp{ x^{ \omega(Z_\alpha) } } }{ \Esp{ x^{ \Peb( \zeta_\Pe(\alpha) ) } } } \times e^{ (x-1)\,\prth{ \log( \alpha - 1) + \zeta_\Pe(\alpha) } } \\
                   & \equivalent{\alpha\to 1} \ 	\Phi_{\widehat{\omega}}(x) \times e^{ (x-1)\,\prth{ \log( \alpha - 1) + \zeta_\Pe(\alpha) } }
\end{align*}

Now, we have (see e.g. \cite[I.1.5, thm. 8 p. 15 \& p. 17-18]{Tenenbaum})
\begin{align}\label{Eq:PrimeZetaEstimate}
\zeta_\Pe( 1 + \varepsilon ) = \log\prth{ \frac{1}{\varepsilon } } - c_0 + o(\varepsilon), 
\qquad \varepsilon > 0, 
\quad c_0 := - \sum_{n \geq 2} \frac{\mu(n)}{n } \log \zeta(n) \approx 0,\!315718
\end{align}

As a result,
\begin{align*}
\frac{ \Esp{ x^{ \omega(Z_\alpha) } } }{ \Esp{ x^{ \Peb( \log\, \prth{  (\alpha - 1)\inv } ) } } }  \tendvers{\alpha }{ 1 } \Phi_{ \widehat{\omega} }(x) \times e^{c_0(x - 1)}
\end{align*}
which itself implies that $ e^{\kappa(x - 1)} \Phi_\omega(x) \Gamma(x) = e^{c_0(x - 1)} \Phi_{ \widehat{\omega} }(x) $. This concludes the proof.
\end{proof}

\medskip

\begin{remark}
To strengthen the analogy with the case of $ C(\sigmab_{\! n}) $ treated in \eqref{Thm:ModPoissonCvTotalNumberCycles}, one can take $ \alpha = \alpha_n = 1 + 1/\log n $. In this case, $  d_{\alpha_n} \stackrel{\Le}{\approx} \exp\prth{ \ee / (\alpha_n - 1) } = n^{ \ee }  $. Such a speed shows that in the transposition of the reasonning of $ C(\sigmab_{\! n}) $, the mod-Poisson speed of convergence must be $ \log \log n + O(1) $ to allow writing $ \log\log(n^\ee) = \log\log n + \log(\ee) $.
\end{remark}

\begin{remark}\label{Rk:BigOvsLittleO:SatheSelberg}
In the same vein as remark~\ref{Rk:BigOvsLittleO}, one can suppose $ o_x(1) $ in place of $ O_x(1/\ln(n)) $ in theorem~\ref{theorem:FacteurGammaSatheSelberg}.
\end{remark}


\begin{remark}[On the convergence \eqref{Thm:SatheSelberg}/\eqref{Thm:SatheSelberg:bis}]
Theorem~\ref{theorem:FacteurGammaSatheSelberg} only provides an ``explanation'' of the appearance of the Gamma factor as coming from the fluctuations of the randomisation $ d_\alpha \sim \deltab\zetab(\alpha) $ when $ \alpha \downarrow 1 $. To turn it into something effective and prove the convergence itself with the desired form of the limiting function, one needs to show that the hypothesis of mod-Poisson convergence with remainder $ O_x(1/\ln n) $ is satisfied (or with remainder $ o_x(1) $ with remark~\ref{Rk:BigOvsLittleO:SatheSelberg}), with an unknown function $ \Phi_\omega $ that is yet to be discovered. Usually, this kind of convergence is proven by making explicit the function, as a result, theorem~\ref{theorem:FacteurGammaSatheSelberg} would only provide a new expression of the limiting function, in the vein of the proof of the Euler formula \eqref{Eq:EulerFormulaGamma}. 

One such case is for instance furnished by Turan in his Doctoral Dissertation (published in Hungarian in 1934) and cited in \cite[\S~0.5 p. 10]{TenenbaumDivisors} and \cite[vol. 2, (2) \& p. 18-19]{Elliott}. Turan shows that for $ x \in \crochet{\frac{1}{\sqrt{2}}, \sqrt{2}} $, one has
\begin{align*}
\Esp{ x^{\omega(U_n) } }  = e^{\ln\ln n \times (x - 1)} \times D(x) (1 + o_x(1))
\end{align*}
with a function $ D $ ``provided by complex integration''. Turan adds that his proof also works for complex $x$ values in the relevant range.

Using a $ o_x(1) $ instead of a $ O_x(1/\ln n) $, one can also use the Halberstam-Richert estimates \cite[ch. 0, thm. 00]{TenenbaumDivisors} to prove the existence of a mod-Poisson limit without computing it. See e.g. \cite[ch. 0.2, thm. 04]{TenenbaumDivisors} for the case of $ \Omega(U_n) $ written under the form~:
\begin{align*}
\Esp{ y^{\Omega(U_n) } }  = \int_0^1 F_{k, y}\prth{ \tfrac{1}{\ln(2nu) }} e^{(y - 1)\ln\ln(2nu)} du 
                              + O_{k, y}\prth{ e^{\ln\ln(2n) \left( \Re(y) - k - 2 \right)}  }
\end{align*}
where $ F_{k, y} $ is a polynomial of degree $k$.
\end{remark}

$ $

\section{The $ \Gamma $-universality class: some examples}\label{Sec:GammaUniversality}

\subsection{\textbf{Goals}}\label{SubSec:GammaUniv:Goals}

Instead of copying the proof of theorem~\ref{theorem:FacteurGammaSatheSelberg} in the numerous cases where the randomised independence structure is proven to exist, we choose instead to collect such cases with the required ingredients (analogues of the Cycle Index Polynomial and Newton's binomial formula interpreted as an independence structure with a randomisation) and to enunciate the theorem without a proof, leaving the task of filling in the details to the interested reader. We will be fast on the reminders, as they can be found in the literature, and will focus on the key properties, i.e.
\begin{enumerate}

\medskip
\item the existence of a ``prime objects decomposition'' in the vein of \eqref{Thm:PrimeDecomposition} or the decomposition of a permutation into products of disjoint cycles, leading to the existence of functionals in the vein of the $k$-cycles $ c_k(\sigma) $ or the $p$-adic valuations $ v_p(n) $ characteristic of the ``prime objects counts'', 

\medskip
\item the existence of an \textit{additive functional} of the prime object counts (e.g. $ C(\sigma) = \sum_k c_k(\sigma) $, $ \omega(n) = \sum_p \Unens{ v_p(n) \geq 1} $ or $ \Omega(n) = \sum_p v_p(n) $, \&c.) that is evaluated in some uniform distribution on the set of objects with a certain size $n$ (e.g. $ U_n \sim \Us(\intcrochet{1, n}) $ or $ \sigmab_{\! n} \sim \Us(\Sg_n) $), and which satisfies a Poisson approximation identity with a particular diverging ``speed'' in the vein of \eqref{Thm:ErdosKacPoisson} ($ C(\sigmab_{\! n}) \approx \Peb(\ln n) $ or $ \omega(U_n) \approx \Peb(\ln\ln n) $),

\medskip
\item the existence of a randomisation identity that expresses in a probabilistic way a ``cycle index polynomial identity'' in the vein of \eqref{Thm:ProbabilisticPolyaCIP} or \eqref{Thm:ProbabilisticArithmeticCIP} with a particular randomisation in the vein of $ g_t $ for \eqref{Thm:ProbabilisticPolyaCIP} and $ d_\alpha $ for \eqref{Thm:ProbabilisticArithmeticCIP},

\medskip
\item the existence of exponential fluctuations for the randomisation when the underlying parameter is set to a limiting value, up to a transformation involving the particular form of the ``speed'', e.g. $ \ln(t\inv) g_t \approx \ee $ or $ \varepsilon \ln d_{1 + \varepsilon} \approx \ee $.

\end{enumerate}

\medskip

These ingredients allow to copy-paste the proof of theorem~\ref{theorem:FacteurGammaSatheSelberg} and show that, up to a factor $ e^{\kappa (x - 1)} $ for a certain $ \kappa $, the limiting mod-Poisson limit, if it exists, is necessarily of the form $ \frac{1}{\Gamma} \times \Phi_{\mathrm{indep.~model}} $.

\medskip
\subsection{Function fields and necklaces}\label{SubSec:GammaUniv:Fq[X]}

\subsubsection{\textbf{Reminders on function fields}}\label{SubSubSec:GammaUniv:Fq[X]:Reminders}

Let $q$ be a power of a prime number. Monic polynomials with coefficients in the field with $ q $ elements $ \Ff_q $ share many similarities with integers and permutations : they have a prime decomposition, the prime polynomials being the irreducible monic polynomials of $ \Ff_q[X] $. One thus has an equivalent of \eqref{Thm:PrimeDecomposition} in the function field setting.

Let $ \Pe\prth{ \Ff_q[X] } $ denote the set of irreducible monic polynomials of $ \Ff_q[X] $. For $ Q \in \Ff_q[X]_m $ and $ P \in \Pe\prth{ \Ff_q[X] } $, denote by $ w_P(Q) $ the $P$(-adic) valuation of $ Q $, namely the exponent of $Q$ in its prime decomposition:
\begin{align}\label{Thm:PrimeDecompositionPolynomials}
Q = \prod_{P \in \Pe( \Ff_q[X] ) } P^{ w_P(Q) }
\end{align}

Denote by $ C_k(Q) $ the number of irreducible factors of degree $k$ of $ Q \in \Ff_q[X]_m $. One has the classical \textit{generalised cyclotomic identity} (see e.g. \cite{MetropolisRota}) valid for $ \abs{t} < 1/q $ and $ \abs{x_k} < 1 $
\begin{align}\label{Thm:CyclotomicIdentity}
\sum_{ Q \in \Ff_q[X]_m } t^{ \deg Q } \prod_{ k \geq 1 } x_k^{C_k(Q) } 
                = \prod_{ P \in \Pe( \Ff_q[X] ) } \frac{1}{1 - x_{\deg P} t^{ \deg P } } 
                = \prod_{k \geq 1 } \prth{ \frac{1}{1 - x_k t^k } }^{ M_q(k) }
\end{align}
where 
\begin{align}\label{Def:NombreDePolynomesIrreductibles}
M_q(k) := \sum_{ P \in \Pe( \Ff_q[X] ) } \Unens{ \deg P = k } = \frac{1}{k} \sum_{ \ell \divise k } \mobius(\ell) q^{ k/\ell }
\end{align}

This formula is the M\"obius inversion of the identity $ \sum_{\ell \divise k} \ell M_q(\ell) = q^k $ that comes from counting in two different ways the roots of some irreducible polynomial in a degree $n$ extension of $ \Ff_q $\footnote{Every element in such an extension is a root of some irreducible polynomial. $ \Ff_q $ itself is the set of roots of the polynomial $ X^q - X $. Here, we consider the polynomial $ X^{q^k} - X $, whose roots are the fixed points of the $k$-th iteration of the Frobenius map $ x \mapsto x^q $. This polynomial factors over $ \Ff_q $ as the product of all irreducible polynomials of degree dividing $k$, hence the formula. }.

There are $ q^n $ monic polynomials of degree $n$ with coefficients in $ \Ff_q $. Let $ g_t \sim \Geom(t) $ and $ Q_n \equiv Q_n^{(q)} \sim \Us\prth{ \ensemble{ Q \in \Ff_q[X]_m  /  \deg Q = n } } $ given by 
\begin{align*}
\Prob{ Q_n = Q } & = \frac{1}{ q^n } \Unens{ \deg Q = n }    
\end{align*}

The original \textit{cyclotomic identity} is obtained by setting $ x_k = (y/t)^k $ in \eqref{Thm:CyclotomicIdentity}. It reads
\begin{align}\label{Thm:CyclotomicIdentitySimple}
\frac{1}{1 - q y} = \prod_{k \geq 1} \prth{ \frac{1}{1 - y^k} }^{ M_q(k) }
\end{align}

Multiplying \eqref{Thm:CyclotomicIdentity} by $ (1 - qt) = \prod_{k \geq 1} (1 - t^k)^{ M_q(k) } $ and replacing $ t $ by $ t/q $, one gets
\begin{align*}
\sum_{n \geq 0 } \Prob{ g_t = n } \Esp{ \prod_{k \geq 1 } x_k^{ C_k(Q_n) } } 
                = \prod_{k \geq 1} \prth{ \Esp{ x_k^{ \Geb((t/q)^k) } } }^{ M_q(k) }
                = \prod_{k \geq 1}  \Esp{ x_k^{ \Neb\Beb( (t/q)^k,\, M_q(k) ) } }    
\end{align*}
where $ \Neb\!\Beb( (t/q)^k, M_q(k) ) \sim \NegBin( (t/q)^k, M_q(k) ) $.

%


One thus gets the equivalent of the CIP theorems \eqref{Thm:ProbabilisticPolyaCIP} and \eqref{Thm:ProbabilisticArithmeticCIP} 
\begin{align}\label{Thm:ProbabilisticPolynomialCIP}
\sum_{n \geq 0 } \Prob{ g_t = n } \Esp{ \prod_{k \geq 1 } x_k^{ C_k(Q_n) } } 
            = \prod_{k \geq 1} \Esp{ x_k^{ \Neb\!\Beb( (t/q)^k,\, M_q(k) )  } } 
\end{align}


\begin{remark}
Using $ \sum_{k \geq 1} k C_k(Q_n) = \deg(Q_n) = n $, one gets the following equality in law analogous to \eqref{Thm:EqLawDimensionGeom}  
\begin{align}\label{Thm:EqLawDegreePol}
\sum_{ k \geq 1} k \Neb\!\Beb( (t/q)^k, M_q(k) ) 
            \eqlaw \sum_{ k \geq 1} k C_k(Q_{ g_t }) 
                 = g_t \sim \Geom(t)
\end{align}
which is the probabilistic translation of \eqref{Thm:CyclotomicIdentitySimple}. This original cyclotomic identity is thus the the analogue of Newton's binomial formula in this context.
\end{remark}

\smallskip
\subsubsection{\textbf{Splitting of $ \Phi_{\omega_q}  $ using the geometric randomisation}}\label{SubSubSec:GammaUniv:Fq[X]:Splitting}

In the same vein as theorem \ref{theorem:FacteurGammaCycles} or \ref{theorem:FacteurGammaSatheSelberg}, one can state the following: 

\begin{shaded}
\begin{theorem}\label{theorem:FacteurGammaPolynomesFq} 
Suppose that $ ( \omega_q(Q_n) )_n $ converges in the mod-Poisson sense at speed $ \log n $ to a certain function $ \Phi_{\omega_q} $ with a speed of convergence given by
\begin{align*}
\frac{ \Esp{ x^{ \omega_q(Q_n) } } }{  \Esp{ x^{\Peb( \log n)} } } = \Phi_{\omega_q}(x) \prth{ 1 + O_x\prth{ \frac{1}{ n} } }
\end{align*}

We define the independent model $ \widehat{\omega}_{q, n} $ by
\begin{align*}
\widehat{\omega}_{q, n} := \sum_{ P \in \Pe( \Ff_q[X] ) , \deg P \leq n } \Beb\prth{ \frac{1}{ \abs{P}_q } } \ \ \ \mbox{ with independent }  \Beb(\abs{P}_q\inv) \sim \Ber_{\{0, 1\}}(\abs{P}_q\inv) 
\end{align*}

Denote by $ \Phi_{ \widehat{\omega}_q }  $ the limiting function of $ (\widehat{\omega}_{q, n})_n $ defined by
\begin{align*}
\Phi_{ \widehat{\omega}_q }(x) := \prod_{P \in \Pe( \Ff_q[X] ) } \prth{ 1 + \frac{x-1}{ \abs{ P }_q } } e^{-\frac{x-1}{ \abs{ P }_q } } 
\end{align*}

Then, necessarily, there exists $ \kappa \in \Rr $ such that $ \Phi_{\omega_q}(x) = \frac{ e^{\kappa (x - 1)} }{\Gamma(x) } \Phi_{ \widehat{\omega}_q }(x)  $.
\end{theorem}
\end{shaded}


\begin{remark}
One can also treat the case of $ \Omega(Q_n) := \sum_{P \in \Pe( \Ff_q[X] ) } w_P(Q_n) $ to get 
\begin{align*}
\frac{ \Esp{ x^{  \Omega_q(Q_n) } } }{ \Esp{ x^{  \Peb(\gamma_{n, q}) } }   }   
\tendvers{ n }{ + \infty }  \Phi_{\Omega_q}\prth{ x } 
                         := \frac{1}{\Gamma(x) } \prod_{ \pi \in \Pe(\Ff_q[X]) }  \frac{1 - \frac{1}{ \abs{\pi}_q } }{  1 - \frac{x }{ \abs{\pi}_q } }  e^{ - \frac{x-1}{ \abs{\pi}_q } }
                         =: \frac{1}{\Gamma(x) }\Phi_{ \widehat{\Omega}_q }(x)
\end{align*}
\end{remark}


\subsubsection{\textbf{Necklaces over an alphabet with $q$ elements}}\label{SubSubSec:GammaUniv:Fq[X]:Other}

The validity of the cyclotomic identity $ (1 - qt)\inv = \prod_{k \geq 1} (1 - t^k)^{ -M_q(k) } $ can extend to $ q \in \Nn^* $ and not necessarily a prime power. The data of the $ (M_q(k))_{k \geq 1} $ that give the number of monic irreducible polynomials when $q$ is a prime power is determined by  \eqref{Def:NombreDePolynomesIrreductibles}, and this expression is a nonnegative integer with which one can try to define an analogous ``prime structure'' amongst $ q^n $ objects.

There exists a natural multiset structure having $ M_q(k) $ different types of objects of weight $k$: the aperiodic words of length $k$ over an alphabet with $q$ elements (see e.g. \cite{MetropolisRota} or \cite[ch. 2.1, ex. 2.12]{ArratiaBarbourTavare} and references cited). The cyclic group acts on such words by rotation. If $ M_q(k) $ designates the number of circular equivalence classes, then, $ k M_q(k) $ is the number of aperiodic words of lenght $k$, hence the relation $ q^n = \sum_{k \divise n} k M_q(k) $. These circular equivalence classes are called \textit{necklaces} or \textit{Lyndon words}.

The analysis in theorem~\ref{theorem:FacteurGammaPolynomesFq} does not suppose anything about $q$ being a prime power and can hence be \textit{de facto} extended to the general case of necklaces, showing the great generality of the randomisation explanation ; this contrasts with \cite{KowalskiNikeghbali2} where the analysis relies on the structure and properties of $ \Ff_q[X] $ and asks the question of generalising the proof to the case of necklaces.


\begin{remark}
When $q$ is the power of a prime, there exists a bijection between necklaces of size $k$ and monic irreducible polynomials of degree $k$, see \cite{Golomb}. This bijection gives the link between the cyclic rotations and the Frobenius map.
\end{remark}


\begin{remark}[Dedekind rings] 

Let $R$ be a Dedekind domain, i.e. for every ideal $I$, the quotient $ R/I $ is finite ; this includes the rings of integers of $ \Zz $ and $ \Ff_q[X] $ or any of their finite extension. Define the norm of an ideal $I$ to be $ N(I) = \abs{R/I} $. Such a norm is multiplicative (i.e. $ N(PQ) = N(P)N(Q) $).
The \textit{Dedekind $ \zeta $ function} of $ R $ is defined by
\begin{align*}
\zeta_R(\alpha) := \sum_{ I \in \operatorname{Ideals}(R) } \frac{1}{ N(I)^\alpha } 
                 = \prod_{ P \in  \operatorname{Spec}(R) } \frac{1}{1 - N(P)^{-\alpha} }
\end{align*}
where $ \operatorname{Spec}(R) $ is the set of prime ideals of $ R $. This Euler product reflects prime factorization of ideals in $ R $. When $ R = \Zz $, one gets the Riemann $ \zeta $ function hence the structure studied in \S~\ref{Sec:Integers}, and when $ R = \Ff_p\crochet{X} $, one gets $ (1 - q\alpha)\inv $, hence \S~\ref{SubSec:GammaUniv:Fq[X]}. This framework thus unifies these two previous cases.

One could also consider the ring of Gauss integers $ \Zz\crochet{i} $ in which the prime ideals are indexed by prime numbers of the form $ 2n + 3 $ and extend the previous framework to this case. 
\end{remark}

\medskip
\subsection{Square free structures : partitions, integers and polynomials over $ \Ff_q $}\label{SubSec:GammaUniv:SquareFree}

\subsubsection{\textbf{Decomposable combinatorial structures}}\label{SubSubSec:GammaUniv:SquareFree:Reminders}

The randomisation in Polya's CIP theorem \eqref{Thm:ProbabilisticPolyaCIP} is governed by the sum $ \sum_{k \geq 1} k C_k $. Indeed, $ \Cb(\sigmab_{\! n}) $ for $ \sigmab_{\! n} \sim \Us(\Sg_n) $ or $ \Cb(Q_n) $ for $ Q_n \sim \Us\prth{ \Ff_q[X]_m \cap \ensemble{ \deg = n } } $ satisfy both the equality $ \sum_{k = 1}^n k C_k = n $ and allow to think the cycle structure or the valuation structure as a conditioning of independent random variables by a linear functional of the $ C_s $'s. 

Other combinatorial structures satisfy such a relation. They are listed in \cite[ch. 2]{ArratiaBarbourTavare}. The case that includes $ \Cb(\sigmab_{\! n}) $ with Poisson random variables is called \textit{assemblies}, the case of $ \Cb(Q_n) $ with negative binomial random variables is included in \textit{multisets}, and a third case of interest concerns \textit{selections}, with binomial random variables. The prototype of this last category of decomposable combinatorial structure is given by strict partitions (namely, partitions of an integers with all parts distincts, see \cite[I-1, ex. 9]{MacDo}), square free integers and square-free polynomials over $ \Ff_q $. We thus call this last category \textit{square free structures}.

\medskip
\subsubsection{\textbf{Square free polynomials over $ \Ff_q $}}\label{SubSubSec:GammaUniv:SquareFree:Fq[X]}

The number of such polynomials is classically known to be $ q^n - q^{n-1} $. This is equivalent to the generating series identity
\begin{align*}
\prod_{n \geq 1} (1 + t^n)^{ M_q(n) } = 1 + \sum_{ n \geq 1 } (q^n - q^{n-1})t^n
\end{align*}
which can be proven using $ 1 + t^n = (1 - t^{2n})/(1 - t^n) $ and the cyclotomic identity \eqref{Thm:CyclotomicIdentitySimple}.

Denote by $ \Phi_q^{(s)} $ the set of square free monic polynomials of $ \Ff_q[X] $. Then, one has the square free cyclotomic identity
\begin{align}\label{Thm:CyclotomicIdentitySquareFree}
\sum_{ Q \in \Phi_q^{(s)} } t^{ \deg Q  } \prod_{k \geq 1} x_k^{C_k(Q) } 
                = \prod_{ P \in \Pe(\Ff_q[X]) } \prth{ 1 + x_{\deg P} t^{ \deg P } } 
                = \prod_{ k \geq 1} \prth{ 1 + x_k t^k }^{ M_q(k) }
\end{align}


Considering $ g_t \sim \Geom(t) $, $ Q_n^{(s)} \sim \Us\prth{ \ensemble{ \deg = n } \cap \Phi_q^{(s)} } $ and
\begin{align*}
\Beb\prth{ M_q(k), \frac{(t/q)^k}{ 1 + (t/q)^k} } \sim \Bin	\prth{ M_q(k), \frac{(t/q)^k}{ 1 + (t/q)^k} }
\end{align*}
one thus gets the equivalent of the CIP theorems \eqref{Thm:ProbabilisticPolyaCIP}, \eqref{Thm:ProbabilisticArithmeticCIP} and \eqref{Thm:ProbabilisticPolynomialCIP}
\begin{align}\label{Thm:ProbabilisticSquareFreePolynomialCIP}
\sum_{n \geq 0 } \Prob{ g_t = n } \Esp{ \prod_{k \geq 1 } x_k^{ C_k\left(Q_n^{(s)}\right) } } 
              = \prod_{k \geq 1} \Esp{ x_k^{ \Beb\left( M_q(k),\, \frac{t^k}{ t^k + q^k} \right) } } 
\end{align}

\medskip
\subsubsection{\textbf{Square free integers}}\label{SubSubSec:GammaUniv:SquareFree:Integers}

An integer $n$ is said to be square free if $ v_p(n) \in \ensemble{0, 1} $ for all $ p \in \Pe $. Denote by $ \Nn_s $ the set of square free integers (starting at $1$).

Let $ (\Beb(w_p))_{p \in \Pe} $ be a sequence of independent Bernoulli random variables of parameters $ w_p \in \crochet{0, 1} $. One can write for all $ x_p \in \crochet{0, 1} $
\begin{align*}
\prod_{p \in \Pe} (1 + (x_p - 1) w_p) 
                             = \Esp{ \prod_{p \in \Pe } x_p^{ \Beb(w_p) } } 
                           & = \sum_{ n \in \Nn^* } \prod_{p \in \Pe} \Prob{ \Beb(w_p) = v_p(n) } \prod_{ p \in \Pe } x_p^{ v_p(n) }  \\
                           & = \sum_{ n \in \Nn_s } \prod_{p \in \Pe} (1 - w_p) \prth{ \frac{w_p}{1 - w_p} }^{ v_p(n) } \prod_{ p \in \Pe } x_p^{ v_p(n) }  \\
                           & =  \prod_{p \in \Pe} (1 - w_p) \times \sum_{ n \in \Nn_s } \prod_{p \in \Pe}  \prth{ \frac{w_p x_p}{1 - w_p} }^{ v_p(n) }  
\end{align*}

Setting $ w_p/(1 - w_p) = p^{-\alpha} $ with $ \alpha > 1 $ and dividing by $ \prod_{p \in \Pe }(1 - w_p) $, one recovers the equivalent of Euler's identity \eqref{Thm:FormuleEulerProduitSommePremiers} for square free integers
\begin{align}\label{Thm:FormuleEulerProduitSommeSquareFree}
\prod_{p \in \Pe} \prth{ 1 + \frac{x_p}{p^\alpha } } = \sum_{ n \in \Nn_s } \frac{1}{n^\alpha} \prod_{p \in \Pe} x_p^{ v_p(n) }  
\end{align}

Define
\begin{align*}
\zeta_{\Nn_s}(\alpha) & := \sum_{ n \in \Nn_s } \frac{1}{n^\alpha } = \prod_{p \in \Pe} \prth{ 1 + \frac{1}{p^\alpha } } \\
\Prob{ Z_\alpha^{(s)} = n }  & := \frac{1}{n^\alpha \zeta_{\Nn_s}(\alpha) } \Unens{n \in \Nn_s}
\end{align*}

This last random variable defines the distribution $ \zetab_{\Nn_s}(\alpha) $. From \eqref{Thm:FormuleEulerProduitSommeSquareFree}, we have $ w_p = \frac{p^{-\alpha}}{1 + p^{-\alpha} } = \frac{1}{p^\alpha + 1} $ and
\begin{align*}
Z_\alpha^{(s)} \eqlaw \prod_{p \in \Pe} p^{ \Beb\left( \frac{1}{p^\alpha + 1} \right) } 
\end{align*}
with independent Bernoulli random variables.


Define $ d_\alpha^{(s)} \sim \deltab\zetab_{\Nn_s}(\alpha) $ by
\begin{align}\label{Def:DeltaZetaDistributionSquareFree}
\Prob{ d_\alpha^{(s)} = n } = \frac{1}{\zeta_{\Nn_s}(\alpha) } \card( \Nn_s \cap \intcrochet{1, n} ) \prth{ \frac{1}{n^\alpha} - \frac{1}{(n + 1)^\alpha} } \Unens{n \in \Nn_s}
\end{align}

It is immediate that if $ d_\alpha^{(s)} $ is independent of $ U_n^{(s)} \sim \Us\prth{\intcrochet{1, n} \cap \Nn_s } $, one gets
\begin{align}\label{Thm:ProbabilisticSquareFreeIntegersCIP}
U_{d_\alpha^{(s)}}^{(s)} \eqlaw Z_\alpha^{(s)}
\end{align}

This last equality in law is thus the CIP identity analogous to \eqref{Thm:ProbabilisticPolyaCIP}, \eqref{Thm:ProbabilisticArithmeticCIP}, \eqref{Thm:ProbabilisticPolynomialCIP} and \eqref{Thm:ProbabilisticSquareFreePolynomialCIP}, as it implies
\begin{align}\label{Thm:ProbabilisticSquareFreeIntegersCIP:bis}
\sum_{ n \geq 1 } \Prob{d_\alpha^{(s)} = n } \Esp{ \prod_{ p \in \Pe } x_p^{ v_p(U_n^{(s)}) } } = \prod_{ p \in \Pe } \Esp{   x_p^{ \Beb\left( \frac{1}{p^\alpha + 1} \right)  } }
\end{align}

\medskip

In the same vein as for $ \deltab\zetab(\alpha) $, one can compute the fluctuations of $\deltab\zetab_{\Nn_s}(\alpha) $:

\medskip

\begin{shaded}
\begin{lemma}[Fluctuations of $ \zetab_{\Nn_s}(\alpha) $ and $ \deltab\zetab_{\Nn_s}(\alpha) $]
Let $ Z_\alpha^{(s)} \sim \zetab_{\Nn_s}(\alpha) $, $ d_\alpha^{(s)} \sim \deltab\zetab_{\Nn_s}(\alpha) $ and $ \ee \sim \Exp(1) $. Then, 
\begin{align}\label{Thm:FluctuationsLoiDeltaZetaSquareFree}
\begin{aligned}
( \alpha - 1 ) \log Z_\alpha^{(s)} \cvlaw{\alpha }{ 1^+ } \ee \\
( \alpha - 1 ) \log d_\alpha^{(s)} \cvlaw{\alpha }{ 1^+ } \ee
\end{aligned}
\end{align}
\end{lemma}
\end{shaded}


\begin{proof}
We have
\begin{align*}
\Prob{ Z_\alpha^{(s)} \leq r }  
                                & =  \frac{1}{\zeta_{\Nn_s}(\alpha)} \sum_{k = 1}^r \frac{1}{k^\alpha}\Unens{ k \in \Nn_s} 
                                  =  \frac{1}{\zeta_{\Nn_s}(\alpha)} \int_1^r \frac{1}{x^\alpha} dM_s(x)
\end{align*}
where $ M_s $ is the density of the set $ \Nn_s $ defined by
\begin{align*}
M_s(x) := \sum_{k \leq x} \Unens{k \in \Nn_s} 
         = \card (\Nn_s \cap \crochet{1, x})
\end{align*}

It is proven in \cite[ch. I.3.7 p. 46]{Tenenbaum} that
\begin{align*}
M_s(x) = \frac{6}{\pi^2} x + O(\sqrt{x})
\end{align*}

As a result, for $ r \gg 1 $, one gets
\begin{align*}
\Prob{ Z_\alpha^{(s)} \leq r }  
                                & =  \frac{6}{\pi^2\zeta_{\Nn_s}(\alpha)} \int_1^r \frac{1}{x^\alpha} dx (1 + o(1)) 
                                = 6\frac{1 - r^{-\alpha + 1}}{\pi^2(\alpha - 1)\zeta_{\Nn_s}(\alpha)}   (1 + o(1))
\end{align*}

Moreover, 
\begin{align*}
\ln \zeta_{\Nn_s}(\alpha) & = \sum_{p \in \Pe} \ln\prth{1 + p^{-\alpha} } 
                    = \sum_{p \in \Pe} \frac{1}{p^\alpha} + \sum_{p \in \Pe} \prth{ \ln\prth{1 + p^{-\alpha} } - p^{-\alpha} } \\
                    & = \zeta_\Pe(\alpha) + c'_0 + o(1), 
                        \qquad \alpha \to 1^+, 
                        \quad c'_0 := \sum_{p \in \Pe} \prth{ \ln\prth{1 + p\inv } - p\inv } \\
                    & = \log\prth{ \frac{1}{\alpha - 1} } + c'_0 - c_0 + o(1) \quad \mbox{ with \eqref{Eq:PrimeZetaEstimate} }
\end{align*}
hence
\begin{align*}
\Prob{ Z_\alpha^{(s)} \leq r } 
                                = \frac{6}{\pi^2 e^{ c'_0 - c_0 } } \prth{1 - r^{-\alpha + 1}}  (1 + o(1)) 
\end{align*}

Setting $ r = \pe{ e^{ t/ (\alpha - 1) } }  $ gives then the result up to the constant $ \frac{6}{\pi^2} e^{-c'_0 + c_0 } $. This constant can only be equal to 1 when setting $ t \to +\infty $.

\medskip

Last, an Abel summation shows that
\begin{align*}
\Prob{ d_\alpha^{(s)} \leq r }  
                                & = \Prob{ Z_\alpha^{(s)} \leq r }  + O\prth{ \frac{1}{\zeta_{\Nn_s}(\alpha) \, (r+1)^\alpha } }
\end{align*}

Setting $ r = \pe{ e^{ t/ (\alpha - 1) } }  $ and using the result for $ Z_\alpha^{(s)} $ gives then the result. 
\end{proof}

\medskip
\subsubsection{\textbf{Splitting of $ \Phi_\Omega $ for square free structures}}\label{SubSubSec:GammaUniv:SquareFree:Splitting}

For square free structures, $ \omega(\cdot) \in \{0, 1\} $ and is hence trivial, but $ \Omega $ is non trival (and is then a sum of independent Bernoulli random variables when evaluated in the distribution on $ \Nn^* $). 

We can state the following theorem, whose proof is a mimic of theorem \ref{theorem:FacteurGammaSatheSelberg} :


\begin{shaded}
\begin{theorem}\label{theorem:FacteurGammaSquareFree}
Let $ \Omega_s(X_n^{(s)}) \in \ensemble{ \Omega_q(Q_n^{(s)}), \Omega(U_n^{(s)}) } $ and $ \gamma_n^{(X)} \in \ensemble{ \gamma_n^{(Q)}, \gamma_n^{(U)} } $ with $ \gamma_n^{(Q)} = \log n $ and $ \gamma_n^{(U)} = \log\log n $. Suppose that $ \Omega_s(X_n^{(s)}) $ converges in the mod-Poisson sense at speed $ \gamma_n^{(X)} $ to a limiting function $ \Phi_{ \Omega_s } $. 


Define the independent models limiting functions by
\begin{align*}
\Phi_{\widehat{\Omega}_s  }(x) := \prod_{a \in \Ae} \Esp{ x^{ \Beb(v_a) } } e^{- (x - 1) v_a}
                                 = \prod_{a \in \Ae} \prth{ 1 + (x - 1) v_a } e^{- (x - 1) v_a}
\end{align*}
with $ \Ae \in \ensemble{\Pe, \Nn^*} $ and $ v_a \in \ensemble{ \frac{1}{a + 1},  \frac{t^a}{1 + t^a} } $.

Then, there exists $ \kappa \in \Rr $ such that
\begin{align*}
\Phi_{ \Omega_s  }(x) = \frac{e^{\kappa(x - 1)} }{ \Gamma(x) } \Phi_{\widehat{\Omega}_s  }(x)
\end{align*}
\end{theorem}
\end{shaded}

\medskip
\subsubsection{\textbf{Note on Strict partitions}}\label{SubSubSec:GammaUniv:SquareFree:NotePartitions}

We did not include this case in the previous theorem, but we mention here the structure for completeness. 

A partition $ \lambda \vdash n $ is said to be strict if all its parts are distinct, namely $ \lambda_1 > \lambda_2 > \dots > \lambda_{\ell(\lambda)} $. Denote by $ \Yy_\infty^{(s)} $ the set of strict partitions, and $ \Yy_n^{(s)} := \Yy_n \cap \Yy_\infty^{(s)} $. The generating function for strict partitions is given by \cite[ch. III.8, p. 252]{MacDo}
\begin{align*}
\sum_{ \lambda \in \Yy_\infty^{(s)} } t^{ \abs{\lambda} } = \prod_{k \geq 1} (1 + t^k)
\end{align*}

There exists a bijection between partitions with odd parts and strict partitions. The first bijection was given by Silvester in \cite{SilvesterBijection} with later simplifications. A bijection given by Lascoux in \cite{LascouxSilvester} allows to find that the following statistics $ \kappa $ is preserved by the bijection:
\begin{align*}
\kappa(\lambda) := \sum_{k \geq 1} (-1)^{k + 1} \lambda_k
\end{align*}

Define for $ \lambda \in \Yy_\infty^{(s)} $
\begin{align*}
\varepsilon_k(\lambda) & := \Unens{ \lambda_{k^*} - \lambda_{ k^* + 1} = 1 } \\
\str(\lambda ) & := \sum_{k \geq 1} \varepsilon_k(\lambda )
\end{align*}
where $ k^* = \Le(k) $ is given by Lascoux's bijection $ \Le $. Note that $ \str(\lambda ) $ is the number of maximal strings of consecutive integers in the word $ \crochet{\lambda_2, \lambda_3, \dots, \lambda_{\ell(\lambda) } } $ (see \cite{LascouxSilvester}). Then, one has the generalised Silvester identity \cite[p. 277]{LascouxSilvester} 
\begin{align*}
\sum_{ \lambda \in \Yy_\infty^{(s)} } t^{ \abs{\lambda} } y^{ \kappa(\lambda) } \prod_{k \geq 1} x_k^{ \varepsilon_k(\lambda) } 
                     = \frac{1}{1 - yt} \prod_{k \geq 1} \prth{ 1 + x_k \frac{ y t^{ 2k + 1 } }{1 - y t^{2 k + 1} } }
\end{align*}

For $ y \in \crochet{0, 1} $, define the measure 
\begin{align}\label{Def:MesureEwensPartitionStricte}
\Pp_{y, n}(\lambda) := \frac{ y^{ \kappa(\lambda) } }{ \crochet{t^n}H\crochet{ \frac{yt}{ 1 - t^2 } }  } \Unens{ \lambda \in \Yy_n^{(s)} }
\end{align}

Note that $ \Pp_{1, n} $ is the uniform measure on $ \Yy_n^{(s)} $. Let $ \lambdab_n^{(s)} \sim \Pp_{y, n} $ and 
\begin{align}\label{Def:RandomisationEwensPartitionStricte}
g_t(y) :\eqlaw \sum_{k \geq 1} (2 k + 1) \Geb(  y t^{2k +1}  )
\end{align}
with independent geometric random variables. This is clear that this last sum converges in $ L^1 $ for $ y t < 1 $.


One thus has the equivalent of the CIP theorems \eqref{Thm:ProbabilisticPolyaCIP}, \eqref{Thm:ProbabilisticArithmeticCIP}, \eqref{Thm:ProbabilisticPolynomialCIP}, \eqref{Thm:ProbabilisticSquareFreePolynomialCIP} and \eqref{Thm:ProbabilisticSquareFreeIntegersCIP}/\eqref{Thm:ProbabilisticSquareFreeIntegersCIP:bis}:
\begin{align}\label{Thm:ProbabilisticStrictPartitionsCIP}
\sum_{n \geq 0} \Prob{g_t(y) = n} \Esp{ \prod_{ k \geq 1 } x_k^{ \varepsilon_k(\lambdab_{n, s} ) } }
                  = \prod_{k \geq 1} \Esp{ x_k^{ \Beb( y t^{ 2k + 1 } ) } }
\end{align}

%
%

\medskip
\subsection{Random matrices of $ GL_n(\Ff_q) $}\label{SubSec:GammaUniv:GLnFq}

\subsubsection{\textbf{Reminders on $ GL_n(\Ff_q) $}}\label{SubSubSec:GammaUniv:GLnFq:Reminders}

If polynomials of $ \Ff_q[X] $ can be thought of as $q$-analogs of integers, it is natural to look for $q$-analogs of permutations to see if theorem \eqref{theorem:FacteurGammaCycles} has a natural $q$-deformed extension in the same vein as theorem \eqref{theorem:FacteurGammaPolynomesFq} is the $q$-deformed version of the Sathe-Selberg theorem \eqref{theorem:FacteurGammaSatheSelberg}.

Such a $q$-deformation of permutations exists and is precisely given by matrices of $ GL_n(\Ff_q) $ when $q$ is a prime power. Loosely speaking, $ \Sg_n \equiv GL_n(\Ff_1) $ or $ GL_n(\Ff_\infty) $, in the sense of a limit on $ q $ towards well-choosen functionals\footnote{
One usually defines formally the ``field with one element''. We consider here that $ \Ff_1\crochet{X} = \ensemble{1} $ to recover the case of $ \Sg_n $.
}. We now remind notations and classical properties of $ GL_n(\Ff_q) $, the main reference being \cite[ch. IV]{MacDo}.

$ GL_n(\Ff_q) $ is a finite group with cardinal given by
\begin{align}\label{Eq:CardinalGLnFq}
\card(GL_n(\Ff_q)) =  \prod_{k = 0}^{n-1} (q^n - q^k)  =  \prod_{k = 1}^n q^{k - 1	} (q^k - 1)
\end{align}

To see the analogy with the cycle decomposition of a permutation of $ \Sg_n $, we look for a characterisation of conjugacy classes of $ GL_n(\Ff_q) $, i.e. equivalence classes for the conjugacy relation~$ \sim $ ($ g_1 \sim g_2 $ iff there exists $ h \in GL_n(\Ff_q) $ s.t. $ g_1 = h\inv g_2 h $). There is a one-to-one correspondence between conjugacy classes in $ M_n(\Ff_q) $ and polypartitions, i.e. families of partitions or partitions-valued functions $ \bbmu :  f \in \Phi_q  \mapsto \bbmu_f \in \Yy  $ that satisfy
\begin{align*}
\norm{\bbmu} := \sum_{f \in \Phi_q} \deg(f) \abs{ \bbmu_f } = n
\end{align*}

We set $ \bbmu_f = (\bbmu_{f, k})_{1 \leq k \leq \ell( \bbmu_f ) } $. This polypartition determines completely a conjugacy class in $ M_n(\Ff_q) $ and a distinguished representant of such a class is given by the block diagonal element $ J(\bbmu) \in M_n(\Ff_q) $ whose blocks are the $ (J(f^{\bbmu_{f, k} } ) )_{f \in \Phi_q, k \geq 1} $. This is the \textit{Jordan-Frobenius decomposition}. Finding the matrix $ J(\bbmu) $ conjugate to a given matrix $ g \in M_n(\Ff_q) $ consists in finding the \textit{rational canonical form} of $g$. We recall that the companion matrix (or Jordan matrix) $ J(f) $ of a monic polynomial $ f := X^n + \sum_{k = 0}^{n - 1} a_k X^k $ is defined by
\begin{align*}
 J(f) := \begin{pmatrix}
   0 &   &          &   & -a_0     \\
   1 & 0 &          &   & -a_1     \\
     & 1 & \ddots   &   & \vdots   \\
     &   & \ddots   & 0 & -a_{n-2} \\
     &   &          & 1 & -a_{n-1}
\end{pmatrix}
\end{align*}

Using the Jordan-Frobenius decomposition, it is easy to write the characteristic (and the minimal) polynomial of a matrix in $ M_n(\Ff_q) $ and see that the matrices that are invertible are those indexed by $ \bbmu : \Phi_q^+ \to \Yy $ (satisfying thus $ \norm{\bbmu}  = n $ and $ \abs{ \bbmu(X) } = 0 $). 

We define $ \bblambda(\alpha) $ as the characteristic polypartition of $ \alpha \in M_n(\Ff_q) $ ; this is the analogue of the cycle type of a permutation, the analogue of the cycles being the Jordan blocks. The total number of such blocks is thus
\begin{align*}
C(\alpha) := \sum_{f \in \Phi_q } \abs{ \bblambda_f(\alpha) }
\end{align*}

For a set of matrices $ S \subset M_n(\Ff_q) $, the cycle index polynomial of $ S $ is the polynomial in variables $ (x_{f, \lambda})_{f \in \Phi_q, \lambda \in \Yy} $ defined by \cite{KungCIP, StongCIP, FulmanCIP}
\begin{align*}
CIP(S) := \frac{1}{\abs{S} } \sum_{ s \in S } \prod_{ f \in \Phi_q } x_{f, \bblambda_f( \alpha ) }
\end{align*}

The following theorem gives the structure of the ``cycle type'' of a matrix in $ GL_n(\Ff_q) $ :


\begin{theorem}[Kung \cite{KungCIP}, Stong \cite{StongCIP}] 
Set $ x_{ 0, \varnothing } := 1 $ and $ CIP(GL_0(\Ff_q)) := 1 $. Let $ \Cent_q(f, \lambda) $ be the centraliser\footnote{
We recall that $ A \oplus B $ designates the diagonal block concatenation of the matrices $ A $ and $B$.
} of $  J(f^{\lambda_1}) \oplus \dots \oplus J(f^{ \lambda_{\ell(\lambda)} }  ) $ in $ GL_{ \abs{\lambda} \deg(f) }(\Ff_q) $ (equal to the the identity matrix if $ \lambda = \varnothing $). One can show that \cite[(2.1)]{BritnellThesis}
\begin{align}\label{Thm:FormuleCardCentralisateur}
\abs{ \Cent_q(f, \lambda) } = q^{ \deg(f) \, \prth{  2 \sum_{k < r} k \, m_k(\lambda) m_r(\lambda) + \sum_k (k-1) m_k(\lambda)^2  }  } \prod_{ k = 1 }^{ \ell(\lambda) } \abs{ GL_{m_k(\lambda)}(\Ff_{ q^{ \deg(f) } }) }
\end{align}

Then, one has the cycle type identity for $ GL_n(\Ff_q) $
\begin{align}\label{Thm:CIPforGLn}
\sum_{n \geq 0 } t^n \, CIP(GL_n(\Ff_q)) = \prod_{f \in \Phi_q^+} \prth{ \sum_{\lambda \in \Yy } \frac{ t^{ \abs{\lambda} \deg(f) } }{ \abs{ \Cent_q(f, \lambda) } } x_{f, \lambda } }
\end{align}
\end{theorem}

We now give a probabilistic interpretation of \eqref{Thm:CIPforGLn}. Setting all variables $ x_{f, \lambda} = 1 $, one gets the identity (for $  \abs{t} < 1 $)
\begin{align*}
\frac{1}{1 - t} = \prod_{f \in \Phi_q^+} \prth{ \sum_{\lambda \in \Yy } \frac{ t^{ \abs{\lambda} \deg(f) } }{ \abs{ \Cent_q(f, \lambda) } }  } =: \prod_{f \in \Phi_q^+} \Ze\prth{ t^{ \deg(f) }, q^{ \deg(f) } }
\end{align*}

Indeed, using \eqref{Thm:FormuleCardCentralisateur}, one sees that $ \abs{ \Cent_q(f, \lambda) } $ is a function of $ q^{\deg(f)} $, hence the choice of notation for $ \Ze\prth{ t^{ \deg(f) }, q^{ \deg(f) } } $. This fact amounts to consider the sole case of $ f = X - 1 $. 

Define the measure
\begin{align}\label{Def:MesureDeSchurSpePrinc}
\Mm_{t, q} := \frac{1}{ \Ze(t, q) } \sum_{ \lambda \in \Yy } \frac{ t^{ \abs{\lambda}   } }{ \abs{ \Cent_q( X - 1, \lambda) } } \delta_\lambda
\end{align}

Let $ U_{n, q} \sim \Us\prth{ GL_n(\Ff_q) } $ and $ \lambdab_{t, q} \sim \Mm_{t, q} $. For the sake of clarity, note $ x_{f, \lambda} =: x_f(\lambda) $ where $ (x_f)_{f \in \Phi_q^+} $ are considered as functions $ x_f : \Yy \to \Rr $. Then, the CIP identity \eqref{Thm:CIPforGLn} becomes
\begin{align}\label{Thm:ProbabilisticGLnCIP}
\sum_{n \geq 0 } \Prob{ g_t = n } \Esp{ \prod_{f \in \Phi_q^+} x_f \prth{ \bblambda_f( U_{n, q} ) } }  = \prod_{f \in \Phi_q^+} \Esp{ x_f\prth{ \lambdab_{ t^{ \deg(f) }, q^{ \deg(f) } } } }
\end{align}
namely, with independent random variables
\begin{align*}
\prth{  \bblambda_f\prth{ U_{g_t, q} }  }_{f \in \Phi_q^+} \eqlaw \prth{ \lambdab_{ t^{ \deg(f) }, q^{ \deg(f) } } }_{f \in \Phi_q^+}
\end{align*}

Moreover, the measures \eqref{Def:MesureDeSchurSpePrinc} are in fact particular \textit{Schur measures}. Indeed, as remarked in \cite[Rk. 2 p. 28]{FulmanGLnIncSubs} or in \cite[ch. IV-6]{MacDo}, using \eqref{Eq:CardinalGLnFq}, one can write the probability $ \Mm_{q, t} $ as 
\begin{align*}
\Prob{ \lambdab_{ t, q} = \lambda } 
              = \frac{ t^{ \abs{ \lambda } } }{ \Ze(t, q) } q^{ 2 n(\lambda) } \prod_{ \square \in \lambda } \prth{ 1 - q^{- h(\square) } }^{-2} 
\end{align*}

Using \cite[ch. I-3 ex. 1-3]{MacDo}, one thus deduces that this last probability is a particular specialisation of a Schur function, namely
\begin{align*}
\Prob{ \lambdab_{ t, q} = \lambda } 
                 = \frac{ t^{ \abs{ \lambda } } }{ \Ze(t, q) } s_\lambda\pleth{ \frac{q\inv}{1 - q\inv } }^2 
                 = \frac{ 1 }{ \Ze(t, q) } s_\lambda\pleth{\sqrt{t}  \frac{q\inv}{1 - q\inv } }^2
\end{align*}

Define
\begin{align*}
\Ae(q) := \frac{q\inv}{1 - q\inv } = \sum_{k \geq 1} q^{-k} = \ensemble{q^{-k}}_{k \geq 1}
\end{align*}

The Cauchy identity gives $ \Ze(t, q) = \sum_\lambda t^{ \abs{ \lambda } } s_\lambda\pleth{\Ae(q)}^2 = H\pleth{ t \Ae(q)\Ae(q) }$ and
\begin{align}\label{Def:MesureDeSchurSpePrincBis}
\Prob{ \lambdab_{ t, q} = \lambda } = \frac{ 1 }{  H\pleth{ t \Ae(q)\Ae(q) } } s_\lambda\pleth{ \sqrt{t}\,\Ae(q) }^2
\end{align}
which is the Schur measure of parameters $ (\sqrt{t}\Ae(q), \sqrt{t}\Ae(q)) $ (see \cite{OkounkovSchurMes}). Remark that the normalisation constant $ \Ze(t, q) $ is in fact a (Cauchy) product, namely
\begin{align*}
H\pleth{t \Ae(q) \Ae(q)   } = \prod_{ k, m \geq 1 } \frac{1}{1 - t q^{ -(k + m) } } = \prod_{ k \geq 1 } \frac{1}{ ( 1 - t q^{ - k } )^k  }
\end{align*}
and in particular, 
\begin{align*}
\abs{ \lambdab_{ t, q} } \eqlaw \sum_{ a, b \in \Ae(q) } \Geb(tab) 
                         = \sum_{k, m \geq 1} \Geb\prth{ t q^{ -(k + m) } }
\end{align*}
with independent geometric random variables, a classical fact for Schur measures (see e.g. \cite[\S~A.2]{BarhoumiMaxIndep}) since, for $ \lambdab_{ \Ae, \Be }  $ a random partition Schur-distributed of parameters $ \Ae $ and $ \Be $, one has
\begin{align*}
\Esp{ x^{ \abs{ \lambdab_{ \Ae, \Be } } } } 
              = \frac{1}{H\pleth{\Ae\Be}} \sum_\lambda x^{\abs{\lambda} } s_\lambda\pleth{\Ae} s_\lambda\pleth{\Be} 
              = H\pleth{(x-1) \Ae\Be } = \prod_{a \in \Ae, b \in \Be} \frac{1 - ab}{1 - x ab}
\end{align*}

Following this computation, we define the (discrete) ``plethystic random variable'' by
\begin{align}\label{Def:PlethysticDiscreteRV}
\Esp{ x^{\Heb\pleth{\Ae}} } := H\pleth{(x - 1) \Ae}  
\end{align}

In particular, $ \Heb\pleth{\Ae \Be} = \sum_{a \in \Ae, b \in \Be } \Geb(ab) $ with independent geometric random variables $ \Geb(ab)$.

\medskip
\subsubsection{\textbf{Splitting of $ \Phi_{C_q} $}}\label{SubSubSec:GammaUniv:GLnFq:Splitting}

In the same vein as theorems \eqref{theorem:FacteurGammaCycles}, 
\eqref{theorem:FacteurGammaSatheSelberg}, 
\eqref{theorem:FacteurGammaPolynomesFq} and 
\eqref{theorem:FacteurGammaSquareFree}, one has the:

\begin{shaded}
\begin{theorem} 
Suppose that $ (C(U_{ n, q }) $ converges in the mod-Poisson sense at speed $ \log n $ to a certain function $ \Phi_{C_q} $. Define an independent model $ \widehat{C}_{q, n} $ as the following sum of independent random variables :
\begin{align*}
\widehat{C}_{q, n} 
         & := \sum_{ f \in \Phi_q^+, \, \deg(f) \leq n   } \Heb\pleth{ \Ae(q^{- \deg(f) } )\Ae(q^{- \deg(f) }) }  
           = \sum_{ f \in \Phi_q^+ , \, \deg(f) \leq n  } \, \sum_{a_f, b_f \in \Ae(q^{- \deg(f) } ) }  \Geb\prth{a_f\, b_f} 
\end{align*}
%
%

The limiting function $ \Phi_{ \widehat{C}_q }  $ of $ (\widehat{C}_{q, n})_n $ is defined by
\begin{align*}
\Phi_{ \widehat{C}_q }(x) :=  \prod_{k \geq 1} \prth{ \frac{1 - q^{ - k } }{ 1 - x q^{ - k }  } e^{ - (x-1) q^{ - k  } } }^{ \!\! N_q(k) }, 
\qquad 
N_q(k) := \sum_{ d \divise k } (M_q(d) - \Unens{ d = 1 } ) \frac{k}{d}
\end{align*}
%
%
%

Then, there exists $ \kappa \in \Rr $ such that $ \Phi_{C_q}(x) = \frac{e^{\kappa (x - 1) } }{\Gamma(x) } \Phi_{ \widehat{C}_q }(x)  $.
\end{theorem}
\end{shaded}

\medskip
\section{Conclusion and perspectives}\label{Sec:Conclusion}

\subsection{Summary}\label{SubSec:Conclusion:Summary}

We sumarise here the different random variables, randomisations and limiting mod-Poisson limiting functions we obtain. The symbol $ \equiv $ is to be understood as ``up to an exponential factor''.
%
%
%

\medskip

\hspace{-0.6cm}
  \begin{tabular}{ | c |  c | c | c | c | c | c |} 
    \hline
    Objects  & 
    Law & 
    Randomisation & 
    CIP & 
    Variable & 
    Speed & 
    Limit \\ 
    \hline
    $ \Sg_n \vphantom{ \Big ( } $ & 
    $ \sigmab_{\! n} \sim \Us(\Sg_n) $ & 
    $ g_t \sim \Geom(t) $ & 
    \eqref{Thm:ProbabilisticPolyaCIP} & 
    $ C(\sigmab_{\! n}) $ & 
    $  \log n $  & 
    $ \Phi_C \equiv \frac{1}{\Gamma }  $ \\ 
    \hline
    $ \intcrochet{1, n} \vphantom{ \Big ( } $ & 
    $ U_n \sim \Us(\intcrochet{1, n}) $ & 
    $ d_\alpha \sim \deltab\zetab(\alpha) $ & 
    \eqref{Thm:ProbabilisticArithmeticCIP}  & 
    $ \omega(U_n) $ & 
    $ \log \log n $  & 
    $ \Phi_\omega \equiv \frac{ \Phi_{\widehat{ \omega }} }{\Gamma } $ \\ 
    \hline
    $ \Ff_q[X] \vphantom{ \Big ( } $ & 
    $ Q_n \sim \Us(\ensemble{ \deg = n}) $ & 
    $ g_t \sim \Geom(t) $ & 
    \eqref{Thm:ProbabilisticPolynomialCIP}  & 
    $ \omega_q(Q_n) $ & 
    $ \log n $  & 
    $ \Phi_{\omega_q} \equiv \frac{ \Phi_{\widehat{\omega}_q}  }{\Gamma }  $ \\ 
    \hline
    $ \Ff_q[X]_{s} \vphantom{ \Big ( } $ & 
    $ Q_n^{(s)} \!\! \sim \Us(\ensemble{ \deg = n }_{\! s}) $ & 
    $ g_t \sim \Geom(t) $ & 
    \eqref{Thm:ProbabilisticSquareFreePolynomialCIP} & 
    $ \Omega_q(Q_n^{(s)}) $ & 
    $ \log n$  & 
    $ \Phi_{\Omega_{q, s}}   \equiv \frac{  \Phi_{\widehat{\Omega}_{q, s}}  }{\Gamma } $ \\ 
    \hline
    $  \Nn_s \vphantom{ \Big ( } $ & 
    $ U_n^{(s)} \sim \Us\prth{\intcrochet{1, n}_s  } $   & 
    $  d_\alpha^{(s)} \sim \deltab\zetab_{\Nn_s}(\alpha) $ &    
    \eqref{Thm:ProbabilisticSquareFreeIntegersCIP:bis}  & 
    $ \Omega(U_n^{(s)}) $ &   
    $ \log\log n $  & 
    $ \Phi_{\Omega_s}  \equiv \frac{\Phi_{ \widehat{\Omega}_s }  }{\Gamma }  $   \\ 
    \hline
  $ GL_n( \Ff_q )  \vphantom{ \Big ( } $ & 
  $ U_{n, q} \sim \Us\prth{ GL_n( \Ff_q ) } $ & 
  $  g_t \sim \Geom(t) $ & 
  \eqref{Thm:ProbabilisticGLnCIP} & 
  $ C( U_{n, q} ) $ &   
  $ \log n $  &   
  $ \Phi_{C_q}   \equiv \frac{ \Phi_{ \widehat{C}_q }   }{\Gamma }  $  \\ 
    \hline
  \end{tabular}

\medskip

$ $

Other structures of interest that were not treated here include: 
colorations of combinatorial structures using wreath products (see e.g. \cite[p. 53]{ArratiaBarbourTavare}) which give hypergeometric distributions, 
Knopfmacher's framework \cite[p. 45]{ArratiaBarbourTavare}, 
statistical mechanics of the Grand Canonical Ensemble \cite{VershikMecaStat, VershikQuantumGas},
other measures on permutations such as the Ewens measure \cite{ArratiaBarbourTavare, FerayMeliotNikeghbali1} or on other structures,
 \&c.

\medskip
\subsection{Perspectives}

We remark that the speed $ \ln n $ is associated with $ g_t $ and the speed $ \ln\ln n $ to the diverse Delta-Zeta distributions, hence an additional logarithmic transformation of the randomisation to get the fluctuations. This asks the question of a super-model interpolating all different distributions, and in particular interpolating the randomisations. Probability theory can give the intuition of the relevant type of objects starting from the relevant type of distribution, if one constructs beforehand such a distribution of course. Such a ``super-model'' would allow to understand not only permutations and polynomials over a finite field in the same framework, but also the integers. Let us remark that the set that interpolate between $ \Sg_n $ and $ \Ff_q[X] $ is, from the probabilistic point of view, $ GL_n(\Ff_q) $. Nevertheless, the laws of the ``cycles'' involve the Schur measure and are thus more complicated than the Poisson or the Negative Binomial distribution. A distribution that would interpolate furthermore with Geometric random variables with parameters involving the primes would also add an additional difficulty.  

Of course, the question remains to find a general paradigm/method that would transpose the computations directly from mod-Poisson converging models to mod-Gaussian or mod-$ \phi $ converging models. While we do not believe that such a paradigm is given by the randomisation feature here presented, we still believe this is a feature that should have its importance and that will find a ``conceptual equivalent'' in the mod-Gaussian framework (but not a straightforward equivalent in the sense of the existence of a randomisation that would provide an independence structure).

\medskip\medskip
\appendix
\section{Properties of the Delta-Zeta distribution}\label{Sec:DeltaZeta}

\subsection{Motivations}\label{SubSec:DeltaZeta:Motivations}

We expressed $ \deltab\zetab(\alpha) $ as a size-biased transform in \S~\ref{SubSec:Integers:DeltaZeta}. We now list some additional properties that allow notably to simulate it \textit{without arithmetical considerations}, using a Markov chain, and to give analogies with the Geometric distribution which is the other randomisation of interest for the appearance of the Gamma factor. Most of these properties are independent randomisations with a Poisson random variables. They are a characterisation of discrete infinitely divisible distributions. For instance, with the notations of \S~\ref{Sec:Notations} and using an independent randomisation,
\begin{align*}
\Neb\!\Beb(\rho, \kappa) & \eqlaw \Peb\prth{ \frac{\rho}{1 - \rho} \, \gammab_\kappa}, 
\qquad \gammab_\kappa \sim \Gammab(\kappa) 
\end{align*}

For $ \kappa = 1 $, one recovers the case $ \Geb(\rho) \eqlaw \Peb\prth{ \frac{\rho}{1 - \rho} \, \ee}$ with $\ee \sim \Exp(1)$. 
From this representation, one can deduce the convergence in distribution
\begin{align*}
\frac{1}{n} \Neb\!\Beb(   e^{-c /n}, \kappa ) 
                 \eqlaw  \frac{1}{n} \Peb\prth{ \tfrac{e^{-c/n}}{1 - e^{-c/n}} \gammab_\kappa }
                 & \cvlaw{n}{+\infty} \frac{1}{c} \gammab_\kappa  
\end{align*} 
using the law of large numbers $ \Peb(\alpha\gamma)/\gamma \to \alpha $ a.s. when $ \gamma\to+\infty $, and the convergence in distribution
\begin{align*}
\Neb\!\Beb\prth{\frac{\lambda}{n}, n} \stackrel{\Le}{\approx}
\Neb\!\Beb\prth{ \frac{1}{1 + \frac{n}{\lambda}}, n} 
                  \eqlaw \Peb\prth{ \frac{ \frac{1}{1 + \frac{n}{\lambda}} }{1 - \frac{1}{1 + \frac{n}{\lambda}} } \, \gammab_n }
                   = \Peb\prth{ \frac{\lambda}{n} \sum_{k = 1}^n \ee_k }
                   \cvlaw{n}{+\infty} \Poisson(\lambda)
\end{align*}
with the law of large numbers $ \gammab_n/n \eqlaw n\inv \sum_{k = 1}^n \ee_k \to 1 $ a.s. when $ n \to +\infty $. One sees that the Negative Binomial distribution interpolates between the Poisson distribution and the Gamma distribution.

\medskip
\subsection{The Zeta and Delta-Zeta distributions as Poisson randomisations}\label{SubSec:DeltaZeta:PoissonRandomisation}


$ $

\begin{shaded}
\begin{theorem}[$ \zetab(\alpha) $ and $ \deltab\zetab(\alpha) $ as randomisations]\label{theorem:Zeta+DeltaZeta=Poisson}
Let $ \gammab_\alpha \sim \Gammab(\alpha) $, $ Z_\alpha \sim \zetab(\alpha) $ and $ \Deb(\alpha) \sim \deltab\zetab(\alpha) $. Define moreover a random variable $ Y_\alpha \in \Rr_+ 	$ with law
\begin{align*}
Y_\alpha \sim \frac{e^{-t}}{1 - e^{- t} } \frac{t^{\alpha - 1} dt}{\Gamma(\alpha)\zeta(\alpha)}
\end{align*}

Then, with independent randomisations and random variables in the RHS
\begin{align}\label{Eq:GeometricRandomisations}
\begin{aligned}
Y_\alpha &\eqlaw \frac{\gammab_\alpha}{Z_\alpha} \\ 
Z_\alpha & \eqlaw 1 + \Geb\prth{e^{-Y_\alpha}} \\
\Deb(\alpha) & \eqlaw 1 + \Neb\!\Beb\prth{e^{-Y_\alpha}, 2 } 
\end{aligned}
\end{align}
and as a result, with independent randomisations and with $ Y_\alpha $ independent of $ \ee $ and $ \gammab_2 $
\begin{align}\label{Eq:PoissonRandomisations}
\begin{aligned} 
Z_\alpha & \eqlaw 1 + \Peb\prth{\frac{e^{-Y_\alpha}}{1 - e^{-Y_\alpha}} \, \ee } \\
\Deb(\alpha) & \eqlaw 1 + \Peb\prth{\frac{e^{-Y_\alpha}}{1 - e^{-Y_\alpha}} \, \gammab_2 } 
\end{aligned}
\end{align}
\end{theorem}
\end{shaded}



\begin{proof}
It is well known that for $ \alpha > 1 $, 
\begin{align*}
\zeta(\alpha) = \int_{\Rr_+} \frac{e^{-t}}{1 - e^{- t} } \frac{t^{\alpha - 1} dt}{\Gamma(\alpha)} 
\end{align*}

As a result, the distribution of $ Y_\alpha $ is well-defined. Its Mellin transform is equal to
\begin{align*}
\Esp{Y_\alpha^s } & = \frac{\Gamma(\alpha + s) \zeta(\alpha + s)}{\Gamma(\alpha) \zeta(\alpha)} \\
               & = \Esp{ \gammab_\alpha^s} \prod_{p \in \Pe} \frac{1 - p^{-\alpha}}{1 - p^{-(\alpha + s) } } \\
               & = \Esp{ \gammab_\alpha^s} \prod_{p \in \Pe} \Esp{ p^{-s \Geb(p^{-\alpha}) } } \\
               & = \Esp{ \gammab_\alpha^s} \Esp{ Z_\alpha^{-s} } 
                 = \Esp{ (\gammab_\alpha Z_\alpha\inv)^s }  
\end{align*}
hence the first result by Fourier-Mellin injectivity. For the second equality in law, one writes
\begin{align*}
\Prob{Z_\alpha = k} & = \frac{1}{\zeta(\alpha)} \frac{1}{k^\alpha}  \\
                   & = \frac{1}{\zeta(\alpha)} \int_{\Rr_+} e^{-kt}  \frac{t^\alpha}{\Gamma(\alpha)} \frac{dt}{t} \\ 
                   & = \frac{1}{\zeta(\alpha)} \int_{\Rr_+} e^{-kt}e^t (1 - e^{-t}) \times e^{-t}(1 - e^{-t})\inv  \frac{t^\alpha}{\Gamma(\alpha)} \frac{dt}{t} \\ 
                   & = \int_{\Rr_+} \Prob{\Geb(e^{-t}) = k - 1 } \Prob{Y_\alpha \in dt} \\
                   & = \Prob{\Geb(e^{-Y_\alpha}) + 1 = k }
\end{align*}
and similarly
\begin{align*}
\Prob{D(\alpha) = k} & = \frac{k}{\zeta(\alpha)} \prth{ \frac{1}{k^\alpha} - \frac{1}{(k +1)^\alpha } } \\
                   & = \frac{k}{\zeta(\alpha)} \int_{\Rr_+} \prth{ e^{-kt} - e^{-(k + 1)t} } \frac{t^\alpha}{\Gamma(\alpha)} \frac{dt}{t} \\
                   & = \frac{k}{\zeta(\alpha)} \int_{\Rr_+} e^{-kt}\prth{ 1 - e^{- t} } \frac{t^\alpha}{\Gamma(\alpha)} \frac{dt}{t} \\
                   & = \int_{\Rr_+} k e^{-kt} e^t\prth{ 1 - e^{- t} }^2 \times e^{-t} \prth{ 1 - e^{- t} }\inv \frac{t^\alpha}{\Gamma(\alpha)\zeta(\alpha)} \frac{dt}{t} \\ 
                   & =: \int_{\Rr_+} \Prob{\Neb\!\Beb(e^{-t}, 2) = k - 1} \Prob{Y_\alpha \in dt} \\
                   & = \Prob{\Neb\!\Beb(e^{-Y_\alpha}, 2) + 1 = k }
\end{align*}
hence the result. The remaining Poissonian equalities in law are consequences of the Poisson representation of the Geometric and Negative Binomial distributions. This concludes the proof.
\end{proof}


\medskip

\begin{remark}
Remark that the second and the first equalities in \eqref{Eq:GeometricRandomisations} imply the fixed point equation in distribution
\begin{align*}
Z_\alpha \eqlaw 1 + \Geb\prth{e^{- \gammab_\alpha/Z_\alpha}}
\end{align*}

As a result, $ Z_\alpha $ is the limit of the following Markov chain (with obvious notations)
\begin{align*}
X_{k + 1} = 1 + \Geb_{k + 1}\prth{e^{- \gammab_\alpha^{(k + 1)}/X_k}}, \qquad X_0 \geq 1
\end{align*}
where $ (\gammab_\alpha^{(k})_k $ and $ (\Geb_k(\cdot))_k $ are i.i.d. sequences that are both independent. Equivalently, using the functional relation $ \Geb(t) \eqlaw \pe{ \frac{\ee}{\ln(t\inv)} } $ and defining an i.i.d. sequence of random variables $ (R_k(\alpha))_{k \geq 1} $ with $ R_1 \eqlaw \frac{\ee}{\gammab_\alpha} $ (with $ \ee $ independent of $ \gammab_\alpha $) one has 
\begin{align*}
X_{k + 1} 
          = 1 + \pe{R_{k + 1}(\alpha) X_k}, \qquad X_0 \geq 1, 
\qquad
Z_\alpha \eqlaw 1 + \pe{ \tfrac{\ee}{\gammab_\alpha} Z_\alpha }
\end{align*}

This simple recursion allows to easily simulate $ Z_\alpha $ \textit{without any arithmetical considerations} (no truncated $ \zeta(\alpha) $ series, etc.) by stopping after a certain number of iterations. It would certainly be interesting to have the speed of convergence of this Markov chain (for instance in total variation distance) to know when to stop iterating, and to see in particular if a cutoff phenomenon is present (see e.g. \cite{DiaconisShahshahaniCutOff}).
\end{remark}

\medskip\medskip
\subsection*{Acknowledgements}

This article is an improved version of some results of the author's Ph.D. thesis. Many thanks are due to my advisor A. Nikeghbali for introducing me to the subject, for constructive criticism and encouragements, to  L. Mutafchiev, J. Najnudel and O. H\'enard for many interesting discussions and corrections, 
and to R. Chhaibi for discussions concerning the Delta-Zeta distribution. Part of this work was done while the author was visiting U. C. Irvine. Many thanks are due to this institution and in particular M. Cranston.

\bibliographystyle{amsplain}

\begin{thebibliography}{10}


\bibitem{ArratiaBarbourTavare}
\rm R. Arratia, A. D. Barbour, and S. Tavaré, 
\it Logarithmic combinatorial structures, a probabilistic approach, 
\rm EMS Monographs in Mathematics, Z\"urich, Europ. Math. Soc. (\textbf{2003}).


\bibitem{BarhoumiModOmega}
\rm Y. Barhoumi-Andr\'eani,
\it A penalised model reproducing the mod-Poisson fluctuations in the Sath\'e-Selberg theorem,
\rm J. Theoret. Probab. 33(2):567-589 (\textbf{2020}).


\bibitem{BarhoumiCUErevisited}
\rm Y. Barhoumi-Andr\'eani,
\it A new approach to the characteristic polynomial of a random unitary matrix,
\rm \href{https://arxiv.org/abs/2011.02465}{arXiv:2011.02465} (\textbf{2020}).


\bibitem{BarhoumiMaxIndep}
\rm Y. Barhoumi-Andr\'eani,
\it Remarkable structures in integrable probability, I: max-independence structures,
\rm \href{https://arxiv.org/abs/2401.08033}{arXiv:2401.08033} (\textbf{2024}).



\bibitem{BarbourKowalskiNikeghbali}
\rm A. D. Barbour, E. Kowalski, A. Nikeghbali,
\it Mod-discrete expansions,
\rm Probab. Theory Related Fields 158(3):859-893 (\textbf{2014}).


\bibitem{BHNY}
\rm P. Bourgade, C. Hughes, A. Nikeghbali, M. Yor,
\it The characteristic polynomial of a random unitary matrix : a probabilistic approach,
\rm Duke Math. J. 145:45-69 (\textbf{2008}).


\bibitem{BritnellThesis}
\rm J. R. Britnell, 
\it Cycle index methods for matrix groups over finite fields, 
\rm Ph.D. thesis, Oxford university (\textbf{2003}).


\bibitem{DelbaenAl}
\rm F. Delbaen, E. Kowalski, A. Nikeghbali,
\it Mod-$ \phi $ convergence,
\rm Int. Math. Res. Not. 2015(11):3445-3485 \url{http://arxiv.org/abs/1107.5657} (\textbf{2015}).


\bibitem{DiaconisShahshahani:Subgroup}
\rm P. Diaconis, M. Shahshahani,
\it The subgroup algorithm for generating uniform random variables, 
\rm Prob. Eng. Inf. Sc., 1:15-32 (\textbf{1987}).


\bibitem{DiaconisShahshahaniCutOff}
\rm P. Diaconis, M. Shahshahani,
\it Generating a random permutation with random transpositions,
\rm Z. Wahrsch. Verw. Gebiete 57:159-179 (\textbf{1981}).


\bibitem{DiaconuGoldfeldHoffstein}
\rm A. Diaconu, D. Goldfeld, J. Hoffstein,
\it Multiple Dirichlet series and moments of Zeta and $L$-functions,
\rm Compos. Math. 139(3):297-360 (\textbf{2003}).


\bibitem{Elliott}
\rm P. D. T. A. Elliott,
\it Probabilistic Number Theory,
\rm vol. I \& II, Springer-Verlag, New York (\textbf{1979}).


\bibitem{ErdosKac}
\rm P. Erd\"os, M. Kac,
\it The Gaussian law of errors in the theory of additive functions,
\rm Proc. N. A. S. 25:206-207 (\textbf{1939}). 


\bibitem{ErdosKac2}
\rm P. Erd\"os, M. Kac,
\it On the Gaussian law of errors in the theory of additive number theoretic functions,
\rm Amer. J. Math. 62:738-742 (\textbf{1940}). 


\bibitem{FerayMeliotNikeghbali1}
\rm V. F\'eray, P.-L. M\'eliot, A. Nikeghbali, 
\it Mod-phi convergence, I: normality zones and precise deviations, 
\rm Springer Briefs in Probability and Mathematical Statistics (\textbf{2016}).


\bibitem{FerayMeliotNikeghbali2}
\rm V. F\'eray, P.-L. M\'eliot, A. Nikeghbali, 
\it Mod-phi convergence, II: Estimates on the speed of convergence, 
\rm in:S\'eminaire de Probabilit\'es L, Lecture Notes in Mathematics, Springer, vol. 2252 (\textbf{2019}).


\bibitem{FulmanCIP}
\rm J. Fulman,
\it A probabilistic approach toward conjugacy classes in the finite general linear and unitary groups,
\rm Journal of Algebra 212:2:557-590 (\textbf{1999}).


\bibitem{FulmanGLnIncSubs}
\rm J. Fulman,
\it $ GL(n, q) $ and increasing subsequences in non-uniform random permutations,
\rm Ann. of Comb. 6:19-32 (\textbf{2002}).


\bibitem{Grandville}
\rm A. Grandville,
\it The anatomy of integers and permutations, 
\rm \url{http://www.dms.umontreal.ca/~ andrew/PDF/Anatomy.pdf} (\textbf{2008}).


\bibitem{Golomb}
\rm S. Golomb,
\it Irreducible polynomials, synchronizing codes, primitive necklaces and cyclotomic algebra,
\rm Proc. Conf. Comb. Math. and its Appl., Univ. of North Carolina Press, 358-370 (\textbf{1969}). 


\bibitem{GonekHughesKeating}
\rm M. Gonek, C. P. Hughes, J. Keating,
\it A hybrid Euler-Hadamard formula for the Riemann Zeta function, 
\rm Duke Math. J. 136:3:507-549 (\textbf{2007}).

\bibitem{Hwang96PGD}
\rm H.-K. Hwang,
\it Large deviations for combinatorial distributions I : central limit theorems, 
\rm Ann. App. Prob., 6:1:297-319 (\textbf{1996}).


\bibitem{Hwang98PGD}
\rm H.-K. Hwang,
\it Large deviations for combinatorial distributions II : local limit theorems, 
\rm Ann. App. Prob., 8:1:163-181 (\textbf{1998}).


\bibitem{Hwang98rate}
\rm H.-K. Hwang,
\it On Convergence Rates in the Central Limit Theorems for Combinatorial Structures, 
\rm Europ. J. of Comb., 19:3:329-343 (\textbf{1998}).


\bibitem{JacodAl}
\rm J. Jacod, E. Kowalski, A. Nikeghbali,
\it Mod-Gaussian convergence: new limit theorems in probability and number theory,
\rm Forum Mathematicum 23:4:835-873 (\textbf{2011}).


\bibitem{KatzSarnak}
\rm N. M. Katz, P. Sarnak, 
\it Random matrices, Frobenius eigenvalues and monodromy, 
\rm AMS Colloquium Publications, vol. 45 (\textbf{1999}).


\bibitem{KowalskiNikeghbali2}
\rm E. Kowalski and A. Nikeghbali,
\it Mod-Poisson convergence in probability and number theory,
\rm Intern. Math. Res. Not. 18:3549-3587 (\textbf{2010}).


\bibitem{KeatingSnaith}
\rm J.P. Keating, N.C. Snaith,
\it Random matrix theory and $ \zeta(1/2 + i t)$,
\rm Comm. Math. Phys. 214:57-89 (\textbf{2000}).


\bibitem{KeatingSnaith2}
\rm J.P. Keating, N.C. Snaith,
\it Random matrix theory and $L$-functions at $s=1/2$,
\rm Comm. Math. Phys. 214:91-110 (\textbf{2000}).


\bibitem{Khintchin}
\rm A. Ya. Khintchin,
\it Limit laws for sums of independent random variables,
\rm ONTI (russian translation), Moscow (\textbf{1938}).


\bibitem{KungCIP}
\rm J. Kung,
\it The cycle structure of a linear transformation over a finite field,
\rm Lin. Alg. Appl. 36:141-155 (\textbf{1981}).


\bibitem{LascouxSilvester}
\rm A. Lascoux,
\it Silvester's bijection between strict and odd partitions,
\rm Discrete mathematics 277:275-278 (\textbf{2004}).


\bibitem{LascouxSym}
\rm A. Lascoux, 
\it Symmetric functions, 
\rm Course Nankai University \url{http://www.emis.de/journals/SLC/wpapers/s68vortrag/ALCoursSf2.pdf} (\textbf{2001}).


\bibitem{MacDo}
\rm I. G. Macdonald, 
\it Symmetric functions and Hall polynomials, 
\rm Oxford Mathematical Monographs, second edition, The Clarendon Press Oxford University Press (\textbf{1995}).


\bibitem{MetropolisRota}
\rm N. Metropolis, G.-C. Rota, 
\it The cyclotomic identity, 
\rm in: Combinatorics and algebra (C. Greene ed.), Contemp. Math. 34:19-27 (\textbf{1983}).


\bibitem{NikeghbaliYorBarnes}
\rm A. Nikeghbali, M. Yor,
\it The Barnes G function and its relations with sums and products of generalized gamma convolution variables, 
\rm Electron. Commun. Probab. 14:396-411 (\textbf{2009}).


\bibitem{OkounkovSchurMes}
\rm A. Okounkov, 
\it Infinite wedge and random partitions, 
\rm Selecta Math. New Ser. 7:57-81 (\textbf{2001}).


\bibitem{Ram1991}
\rm A. Ram, 
\it A Frobenius formula for the characters of the Hecke algebra, 
\rm Invent. math. 106:61-488 (\textbf{1991}).


\bibitem{RenyiTuran}
\rm A. Rényi, P. T\'uran ,
\it On a theorem of Erd\"os-Kac,
\rm Acta Arithmetica 4:71-84 (\textbf{1958}).


\bibitem{Sathe}
\rm L. G. Sathe, 
\it On a problem of Hardy on the distribution of integers having a given number of prime factors, 
\rm J. Indian Math. Soc. 17(1):63-141 (\textbf{1953}), 18(1):27-81 (\textbf{1954}).


\bibitem{SelbergSathe}
\rm A. Selberg,
\it Note on a paper by L.G. Sathe,
\rm J. Indian Math. Soc. 18(1):83-87 (\textbf{1954}). 


\bibitem{SheppLloyd1966}
\rm L. A. Shepp and S. P. Lloyd, 
\it Ordered cycle lengths in a random permutation, 
\rm \TAMS 121:340-357 (\textbf{1966}).


\bibitem{SilvesterBijection}
\rm J. J. Sylvester,
\it A constructive theory of partitions, arranged in three acts, an interact and exodion,
\rm Amer. J. Math. 5:251-330 (\textbf{1882}) [reprinted in: The Collected Papers 4:1-83, Cambridge University Press (\textbf{1912}) ; reprinted by Chelsea, New York, (\textbf{1974})].


\bibitem{StongCIP}
\rm R. Stong,
\it Some asymptotic results on finite vector spaces,
\rm Adv. Appl. Math. 9:167-199 (\textbf{1988}). 


\bibitem{TaoBlog}
\rm T. Tao,
\it Selberg's limit theorem for the Riemann zeta function on the critical line,
\rm Structure and randomness, t. 2, \url{http://terrytao.wordpress.com/2009/07/12/selbergs-limit-theorem} \url{-for-the-riemann-zeta-function-on-the-critical-line}. 


\bibitem{Tenenbaum}
\rm G. Tenenbaum,
\it Introduction to analytic and probabilistic number theory,
\rm Cambridge Studies in Advanced Mathematics, Cambridge University Press (\textbf{1995}). 


\bibitem{TenenbaumDivisors}
\rm G. Tenenbaum,
\it Divisors,
\rm Cambridge tracts in mathematics, vol. 90, Cambridge University Press (\textbf{1988}).


\bibitem{VershikMecaStat} 
\rm A. M. Vershik,
\it Statistical mechanics of combinatorial partitions, and their limit configurations, 
\rm Funktsional. Anal. i Prilozhen. 30:2:19-39 (\textbf{1996}) [English translation: Funct. Anal. Appl. 30:2:90-105 (\textbf{1996})].


\bibitem{VershikQuantumGas}
\rm A. M. Vershik,
\it Limit distribution of the energy of a quantum ideal gas from the point of view of the theory of partitions of natural numbers, 
\rm Uspekhi Mat. Nauk 30:2:139-146 (\textbf{1997}) [English translation: Russian Math. Surveys 52:2:379-386 (\textbf{1997})]


\bibitem{WhittakerWatson}
\rm E.T. Whittaker, G.N. Watson,
\it A course in modern analysis,
\rm 4th Edition, Cambridge Math. Library, Cambridge University Press (\textbf{1996}). 

$ $

\end{thebibliography}


\end{document}